\documentclass[reqno,11pt]{amsart}

\usepackage{amsmath,amsthm,amssymb,amsfonts,epsfig,color,graphicx,enumerate,accents,subfigure,mathtools}
\usepackage{hyperref}
\usepackage[a4paper,margin=3truecm]{geometry}
\usepackage{aliascnt}
\usepackage{cleveref}
\usepackage{todonotes}

\newtheorem{thr}{Theorem}[section]
\newaliascnt{prop}{thr}
\newtheorem{prop}[prop]{Proposition}
\aliascntresetthe{prop}
\newaliascnt{lem}{thr}
\newtheorem{lem}[lem]{Lemma}
\aliascntresetthe{lem}
\newaliascnt{lemma}{thr}
\newtheorem{lemma}[lemma]{Lemma}
\aliascntresetthe{lemma}
\newaliascnt{cor}{thr}
\newtheorem{cor}[cor]{Corollary}
\aliascntresetthe{cor}

\theoremstyle{remark}
\newaliascnt{rmk}{thr}
\newtheorem{rmk}[rmk]{Remark}
\aliascntresetthe{rmk}

\theoremstyle{plain}
\newaliascnt{assumption}{thr}
\newtheorem{assumption}[assumption]{Assumption}
\aliascntresetthe{assumption}

\theoremstyle{definition}
\newaliascnt{df}{thr}
\newtheorem{df}[df]{Definition}
\aliascntresetthe{df}
\newaliascnt{exe}{thr}
\newtheorem{exe}[exe]{Example}
\aliascntresetthe{exe}
\theoremstyle{plain}
\newtheorem*{ass*}{\textbf{\textup{(A)}}}
\newtheorem{thmx}{Theorem}

\newcommand{\supp}{\textup{supp}}

\newcommand{\I}{\mathcal{I}}
\newcommand{\J}{\mathcal{J}}
\newcommand{\C}{\textup{C}}

\newcommand{\D}{\textup{D}}

\newcommand{\M}{\mathcal{M}}

\renewcommand{\P}{\mathcal{P}}
\renewcommand{\d}{\, \textup{d}}

\newcommand{\R}{\mathbb{R}}

\newcommand{\N}{\mathbb{N}}
\newcommand{\weakto}{\rightharpoonup}

\renewcommand{\span}{\operatorname{span}}

\newif\ifshowcorr
\showcorrfalse
\newcommand{\corr}[1]{\ifshowcorr\textcolor{blue}{#1}\else#1\fi}

\title[Moment approximations in optimal transport]{Moment approximations in optimal transport: geometric estimates and a semiclassical limit}

\author{Luigi De Pascale}
\address{Dipartimento di Matematica e Informatica, Università di Firenze \\ Viale Morgagni 67/a, 50134 Firenze \\ Italy}
\email{luigi.depascale@unifi.it}

\author{Luca Nenna}
\address{Universit\'e Paris-Saclay, CNRS, Laboratoire de math\'ematiques d'Orsay, ParMA, Inria Saclay, 91405, Orsay, France. Institut Universitaire de France (IUF). }
\email{luca.nenna@universite-paris-saclay.fr}

\author{Igor Pinheiro}
\address{Dipartimento di Matematica e Informatica ``Ulisse Dini" \\ Università degli Studi di Firenze}
\email{igorvinicius.pereirapinheiro@unifi.it}

\date{\today}

\crefname{thr}{theorem}{theorems}
\Crefname{thr}{Theorem}{Theorems}
\crefname{prop}{proposition}{propositions}
\Crefname{prop}{Proposition}{Propositions}
\crefname{lem}{lemma}{lemmas}
\Crefname{lem}{Lemma}{Lemmas}
\crefname{lemma}{lemma}{lemmas}
\Crefname{lemma}{Lemma}{Lemmas}
\crefname{cor}{corollary}{corollaries}
\Crefname{cor}{Corollary}{Corollaries}
\crefname{rmk}{remark}{remarks}
\Crefname{rmk}{Remark}{Remarks}
\crefname{assumption}{assumption}{assumptions}
\Crefname{assumption}{Assumption}{Assumptions}
\crefname{df}{definition}{definitions}
\Crefname{df}{Definition}{Definitions}
\crefname{exe}{example}{examples}
\Crefname{exe}{Example}{Examples}
\crefname{thmx}{theorem}{theorems}
\Crefname{thmx}{Theorem}{Theorems}

\begin{document}

\begin{abstract}
We study finite moment approximations of linear optimization problems over measures. Under suitable assumptions, we establish the existence of discrete minimizers with a controlled number of support points, duality and $\Gamma$-convergence as the moment approximation is refined.\corr{ This framework accommodates several transport problems, in which marginal constraints are replaced by finitely many moment constraints}. For moment functions localized at scale $r_n$, we prove an estimate controlling $W_\infty$ by $W_p$ up to an additive error of order $r_n$, which also applies to discrete measures. For repulsive costs, we show that minimizers remain separated from the diagonal, uniformly in the number of moments, and we derive consequences for dual potentials. Finally, under a moment-correction assumption, we prove that the fermionic minimum energies converge in the semiclassical limit to the moment-constrained Coulomb transport value for a fixed family of moment functions.
\end{abstract}


\maketitle

\section{Introduction}\label{sec:introduction}

The numerical approximation of optimal transport problems becomes particularly difficult when several marginals are prescribed. A common approach consists in approximating each marginal by a discrete measure and solving the resulting finite-dimensional problem. If each of $M$ marginals is supported on $n$ points, the transport plan has up to $n^M$ possible configurations. This difficulty motivates approximations that allow the support of the transport plan to be sparse \cite{AlfCoyEhrLom2021Approx,alfonsi2022constrained}.

In \cite{AlfCoyEhrLom2021Approx}, the authors study one such approximation, which consists in replacing the marginal constraints by finitely many moment constraints. The prescribed marginals enter the problem through their integrals against a chosen family of test functions, called moment functions. Under suitable assumptions, the relaxed problem admits a discrete minimizer whose number of support points is controlled by the number of constraints. Both the positions of these points and their weights can then be treated as variables in the approximation.

The choice of moment functions is consequently an essential part of the approximation. The relaxed constraints allow measures whose marginals differ from the prescribed ones, and it is necessary to understand what information about those marginals is retained. In this paper, we study finite moment approximations of linear optimization problems over measures within a common framework, and we investigate which geometric information on the marginals is retained by the relaxation. To fix ideas, we first describe the approximation in the Monge--Kantorovich case.

Let $X$ and $Y$ be Polish spaces, and let $\mu\in\P(X)$ and $\nu\in\P(Y)$. Given a lower semicontinuous cost $c:X\times Y\to\R\cup\{+\infty\}$, bounded from below, the Monge--Kantorovich problem consists in minimizing
\[
    \C(\pi)=\int_{X\times Y}c(x,y)\,\d\pi(x,y)
\]
among probability measures $\pi$ whose first and second marginals are $\mu$ and $\nu$, respectively. We denote the set of all such measures by $\Pi(\mu,\nu)$.
Equivalently, these marginal constraints can be expressed using test functions:
\[
\begin{aligned}
    \int_{X\times Y}\varphi(x)\,\d\pi(x,y)
    &=\int_X\varphi\,\d\mu
    &&\text{for every }\varphi\in C_b(X),\\
    \int_{X\times Y}\psi(y)\,\d\pi(x,y)
    &=\int_Y\psi\,\d\nu
    &&\text{for every }\psi\in C_b(Y).
\end{aligned}
\]
The moment relaxation consists in requiring these identities only for selected finite families of test functions. For each $n\in\N$, let $\I^n\subset C_b(X)$ and $\J^n\subset C_b(Y)$ be finite families of moment functions. We denote by $\Pi^n(\mu,\nu)$ the set of probability measures satisfying the preceding identities for every $\varphi\in\I^n$ and $\psi\in\J^n$. The moment-constrained Monge--Kantorovich problem is then
\begin{equation}\label{eq:mcot-problem}
    \inf_{\pi\in\Pi^n(\mu,\nu)}\C(\pi).
\end{equation}
Since $\Pi(\mu,\nu)\subseteq\Pi^n(\mu,\nu)$, the infimum of the relaxed problem is no larger than that of the original problem; when the moment functions are refined so as to fully characterize $\mu$ and $\nu$, the relaxed values converge to the original one (see \Cref{thm:ot-convergence}).

Even when the relaxation preserves the optimal value, it may admit minimizers with very different supports, as the following example shows.

\begin{exe}\label{exe:polynomial-geometry}
Let $X=Y=[0,1]$, consider the attractive cost $c(x,y)=|x-y|^2$, and take
\[
    \mu=\nu=\frac12\delta_{1/4}+\frac12\delta_{3/4}.
\]
The unique minimizer of the original problem is
\[
    \pi^\star=\frac12\delta_{(1/4,1/4)}+\frac12\delta_{(3/4,3/4)}.
\]
For $\I^3=\J^3=\{x\mapsto x^k:k=1,2,3\}$, a minimizer of the moment-constrained problem is
\[
    \pi^3=\frac18\delta_{(0,0)}+\frac34\delta_{(1/2,1/2)}+\frac18\delta_{(1,1)},
    \qquad \C(\pi^3)=0.
\]
Both its marginals are
\[
    \frac18\delta_0+\frac34\delta_{1/2}+\frac18\delta_1,
\]
which have the same first three moments as $\mu$. Nevertheless, their support is disjoint from that of $\mu$, and $\supp(\pi^3)\cap\supp(\pi^\star)=\varnothing$.
\end{exe}

\Cref{exe:polynomial-geometry} shows that the optimal value alone does not detect how far the relaxed marginals are from the prescribed ones: the moment functions have to be chosen so as to retain geometric information. We now present our main results in a simplified form, together with their consequences; precise statements and assumptions are given in the subsequent sections. Throughout, for a cost $c$ we write $\C(\pi)=\int c\,\d\pi$.

\subsection*{Geometric estimates}

We say that the moment functions $\varphi_1,\dots,\varphi_n$ are \emph{localized at scale $r_n$} on the support of $\mu$ if they form a partition of unity on $\supp(\mu)$ made of functions supported in balls of radius $r_n$. Any probability measure having the same integrals as $\mu$ against $\varphi_1,\dots,\varphi_n$ is called a \emph{relaxed marginal} of $\mu$.

\begin{thmx}[Geometric estimates]\label{thm:intro-geometric}
Assume that the moment functions are localized at scale $r_n$. Then every relaxed marginal $\mu_n$ of $\mu$ satisfies $W_\infty(\mu,\mu_n)\leq2r_n$. Moreover, if $\mu$ has a density bounded from below on a bounded Lipschitz domain $\Omega\subset\R^d$ and $\nu$ is concentrated on $\overline\Omega$, then, for all relaxed marginals $\mu_n$ of $\mu$ and $\nu_n$ of $\nu$,
\[
    W_\infty(\mu_n,\nu_n)\leq 4r_n+C\,W_p(\mu_n,\nu_n)^{\frac{p}{p+d}}.
\]
\end{thmx}

The first estimate rules out the behavior of \Cref{exe:polynomial-geometry}: with localized moment functions, the relaxation cannot move any portion of mass farther than $2r_n$, so that the relaxed marginals converge to the prescribed ones in the strongest Wasserstein sense as $r_n\to0$, and not only narrowly. The second estimate compares an averaged distance with a uniform one. If $\mu_n$ and $\nu_n$ are close on average, namely $W_p(\mu_n,\nu_n)$ is small, then they can be coupled without moving any point farther than $C\,W_p(\mu_n,\nu_n)^{p/(p+d)}$, up to the resolution $4r_n$ of the moment functions. Estimates of this type are classical when one of the measures has a density bounded from below \cite{BoucJimeRaje2007PAMS,JylhaRajala2016}. Here, the constants only depend on the prescribed marginal $\mu$, so that the estimate applies to the relaxed marginals, which are typically discrete. See \Cref{lem:moment-localization}, \Cref{thm:wp-winf-bound} and \Cref{cor:wp-winf-density} for the precise statements.

\subsection*{Repulsive costs}
An important motivation for these sparse discretizations comes from multi-marginal optimal transport problems with Coulomb cost, which arise in the strictly correlated electrons limit of density functional theory \cite{ButDepGorPRA2012,CFK13} and for which moment-constrained approximations were introduced in \cite{AlfCoyEhrLom2021Approx,alfonsi2022constrained,ehrlacher2026sparse}. For such costs, it is natural to ask whether the support of a minimizer of the moment relaxation remains at a positive distance from the diagonal, where the cost is singular.

We consider $M\geq2$ marginals all equal to $\rho\in\P(X)$ and costs of the form $c(x_1,\dots,x_M)=\sum_{i<j}h(d(x_i,x_j))$, where $h$ is decreasing and blows up at $0$, such as the Coulomb cost $h(r)=1/r$ on $\R^3$. The moment relaxation only prescribes the integrals of the marginals against finitely many moment functions, which need not be localized.

\begin{thmx}[Separation from the diagonal]\label{thm:intro-diagonal}
Assume that $\rho$ does not concentrate too much mass on small balls. Then there exists $\alpha^\star>0$, independent of the number $n$ of moment functions, such that every minimizer of the moment-constrained problem gives no mass to configurations in which two points are closer than $\alpha^\star$.
\end{thmx}

 This has two consequences. First, since minimizers never see the singularity, the cost can be replaced by a bounded continuous truncation without changing the minimizers (\Cref{thm:repulsive-truncation}); this justifies the use of bounded costs in numerical computations. Second, when $X$ is compact and the moment functions are localized, the dual problem admits a symmetric maximizer in the span of the moment functions, with explicit bounds on its coefficients (\Cref{thm:repulsive-potentials}), in the spirit of the bounds on Kantorovich potentials for the Coulomb cost \cite{ColoDiMaStra2019SIMA,ButtChamDePa2018AMO}. See \Cref{thm:diagonal-estimate} and \Cref{cor:diagonal-uniform} for the precise statements.

\subsection*{A semiclassical limit}

Our third main result connects the moment-constrained Coulomb problem with its quantum counterpart. In density functional theory, the multi-marginal problem with Coulomb cost arises as the semiclassical limit of the Lieb functional \cite{CFK13,BindDePa2017JEP,Lewi2018CRM,CotarFriKlu2018}: when the kinetic energy is multiplied by a small parameter $\varepsilon$, the minimal energy of $M$ fermions with a prescribed density converges to the optimal transport value. We consider the same quantum problem, where the density constraint is replaced by finitely many moment constraints on the $M$-particle law, as in \cite{ehrlacher2026sparse}.

\begin{thmx}[Semiclassical limit]\label{thm:intro-semiclassical}
Fix $n\in\N$ and let $c$ be the Coulomb cost. Under a technical assumption ensuring that the moment constraints can be corrected, the minimal energy
\[
    \inf_{\Gamma}\left\{
        \frac{\varepsilon}{2}\operatorname{Tr}(-\Delta\Gamma)+\C(\pi_\Gamma)
    \right\},
\]
where the infimum is taken over the fermionic states $\Gamma$ whose $M$-particle law $\pi_\Gamma$ satisfies the moment constraints, converges, as $\varepsilon\to0$, to the minimum of $\C$ under the same moment constraints. Moreover, the laws of almost minimizing states converge, up to subsequences, to minimizers of the moment-constrained transport problem.
\end{thmx}

In other words, imposing finitely many moments and taking the semiclassical limit are compatible operations: the moment-constrained transport problem studied in the previous results is exactly the strictly correlated limit of the moment-constrained quantum problem. This justifies the use of moment-constrained optimal transport as an approximation of the quantum problem in the strongly correlated regime. The result follows from a $\Gamma$-convergence statement (\Cref{thm:semiclassical-gamma}). The construction of recovery sequences relies on three steps: we approximate a transport plan by finitely supported plans with finite cost, which are therefore separated from the diagonal; we replace each atom by a smooth fermionic wave packet, which keeps the Coulomb energy under control; finally, we restore the moment constraints exactly through convex combinations with finitely many fixed states, which is where the technical assumption (\Cref{ass:moment-correctability}) enters.

\subsection*{Related works}

Finite moment relaxations of the marginal constraints were introduced in \cite{AlfCoyEhrLom2021Approx} for multi-marginal optimal transport, where the existence of sparse minimizers and the convergence as the number of moments increases are proved, and numerical methods are proposed; see also \cite{alfonsi2022constrained} for symmetric multi-marginal problems and \cite{ehrlacher2026sparse} for the corresponding approximation of the Lieb functional. Other forms of moment constraints in optimal transport have been considered more recently. Additional linear constraints on the transport plan, with exact marginals, are studied in \cite{Zae2015Monge}. In \cite{carlier2025weak}, weak optimal transport problems with hard and soft moment constraints on the disintegration of the plan, which include martingale transport, are studied, and a Slater-type qualification condition yielding dual attainment is introduced. In \cite{kolesnikov2026duality}, a martingale transport problem in which the second marginal is only prescribed through finitely many moments is considered, and a finite-dimensional dual problem without duality gap is derived under an interiority condition on the moment vector, which plays a role similar to our full-support assumption in \Cref{thm:duality}. Finally, the extreme points of moment sets \cite{Wink1988MOR,Pine2016MMOR} and the Richter--Tchakaloff theorem \cite{BayerTeichmann2006} are classical tools of the theory of moments \cite{schmudgen2017moment,DiDioSchm2022JMAA}.

Compared with these works, \Cref{sec:moment-theory} gathers in a single framework the existence of minimizers, their discreteness and $\Gamma$-convergence for a linear problem over nonnegative measures with finitely many moment inequalities and a confinement constraint, which are then applied to all the problems of \Cref{app:other-applications}. As concerns the duality result, \Cref{thm:duality}, which holds under a full-support condition and reduces the dual problem to a finite-dimensional maximization over the coefficients of the moment functions (see \Cref{rmk:finite-dim-dual}), is, to the best of our knowledge, new in this generality. 
\subsection*{General framework and organization of the paper}
Several variational problems share the linear moment structure used above. In \Cref{sec:moment-theory}, we therefore study a general linear optimization problem over nonnegative measures of fixed mass, subject to moment equalities and inequalities. This formulation allows the basic approximation results to be proved once and then applied through the choice of the underlying space, cost and moment functions. Under suitable assumptions, we establish existence of minimizers, obtain discrete minimizers through the Richter--Tchakaloff theorem and extreme-point arguments, and derive duality and optimality conditions. We also prove $\Gamma$-convergence as the moment approximation is refined. The section concludes with the Monge--Kantorovich case, including convergence of the optimal values and subsequential convergence of minimizers when the moment functions fully characterize the prescribed marginals.

In \Cref{sec:wp-winf}, we study how localized moment functions control the geometry of the relaxed marginals and prove the estimate relating $W_p$ and $W_\infty$. \Cref{sec:repulsive-costs} establishes the diagonal exclusion estimates and their consequences for cost truncation and dual potentials. \Cref{sec:semiclassical-gamma} discusses the semiclassical application through the smoothing of atomic measures and the correction of their moments.

\Cref{app:other-applications} shows how the general theory applies to partial optimal transport, JKO steps with linear potential energies, martingale optimal transport on compact sets, capacity-constrained optimal transport, and Wasserstein-on-Wasserstein transport. For each problem, we identify the finite moment relaxation and state the corresponding existence, duality and convergence results. These applications explain the scope of the abstract formulation: moment inequalities describe partial transport and capacity constraints, additional equalities encode the martingale condition, and the underlying space may itself be a space of probability measures.

\section{General theory of moment-constrained optimization}
\label{sec:moment-theory}

In this section, we study the abstract moment-constrained problem underlying all the applications of this paper; the Monge--Kantorovich case is recovered in \Cref{sec:mk}.\corr{ For background on the general moment problem, we refer to \cite{schmudgen2017moment,DiDioSchm2022JMAA}.}

\subsection{Formulation}\label{sec:formulation}

Let $Z$ be a Polish space, let $\M_+(Z)$ denote the set of finite nonnegative Borel measures on $Z$, and let
$c:Z\to \mathbb{R} \cup \{+\infty\}$ be a lower semicontinuous cost function bounded from below.
Consider the functional
\[
\C:\mathcal M_+(Z)\to \mathbb{R}\cup\{+\infty\},
\qquad
\C(\pi)=\int_Z c\,d\pi.
\]
Let $w_1,\dots, w_n \in C_b(Z)$  be bounded continuous functions. For a fixed, 
coercive\footnote{Every sublevel set of $\Theta$, $\{\Theta \leq R\}$, is compact for every $R\geq0$.}, lower semicontinuous function $\Theta:Z\to [0,+\infty]$, and for
$A>0$ and $m> 0$ we define the following variational problem\footnote{If $\Pi^n_A = \emptyset$, we assume the infimum to be $+\infty$.}
\begin{equation}\label{eq:moment-problem}
\inf_{\pi\in \Pi_A^n} \C(\pi).
\end{equation}where
\begin{equation}\label{eq:moment-set}
\begin{aligned}
\Pi_A^n \coloneqq \{\, &
\pi\in\mathcal M_+(Z):\ \pi(Z)=m,\\
&\int_Z w_i\,d\pi\geq0\quad(i=1,\dots,n),
\quad \int_Z\Theta\,d\pi\leq A\,\}.
\end{aligned}
\end{equation}

\subsection{Existence}

\begin{thr}\label{thm:existence}
     If the infimum in \eqref{eq:moment-problem} is finite, then there exists at least one minimizer. 
\end{thr}
\begin{proof}
As \eqref{eq:moment-problem} is finite we have that $\Pi^n_A \neq \emptyset$. Let $(\pi_k)_{k\in \N} \subset \Pi_A^n $ be a minimizing sequence. Since \[\int \Theta \d \pi_k \leq A \text{ for all } k\in \N\]  and $\Theta$ is coercive, the sequence is tight. Passing to a subsequence without changing notation we have $\pi_k \weakto \pi_\infty$ and $\pi_\infty \in \mathcal{M}_+(Z)$. Let us show that $\pi_\infty$ is admissible. From weak convergence, $\pi_\infty(Z) = m$ and from the lower semicontinuity of $\Theta$, $\int \Theta d\pi_\infty \leq A$. Moreover, the functional $\pi \mapsto \int w \d \pi$ is continuous, as $w\in C_b(Z)$, therefore weak convergence again guarantees that $\int_{Z} w_i \d\pi_\infty \geq 0$ for every $1\leq i \leq n$ and thus 
\[\int c \d\pi_\infty \geq \inf_{\pi \in \Pi^n_A} \int_Z c\d \pi .\] The other inequality follows from the lower semicontinuity of the cost. The same argument, applied to an arbitrary sequence in $\Pi^n_A$, shows that $\Pi^n_A$ is narrowly compact.
\end{proof}

\subsection{Characterization}\label{sec:structure}

\begin{thr}[Richter--Tchakaloff \cite{BayerTeichmann2006}]\label{thm:tchakaloff} Let $Z$ be a Polish space and let $\pi$ be a nonzero finite nonnegative Borel measure concentrated on a Borel set $E\subset Z$. Let $N \in \N$ and let $\Lambda: Z \to \R^{N}$ be a Borel measurable map. Assume that \begin{equation}
    \int_{\R^N} \| u \| \d (\Lambda_{\#}\pi)(u) = \int_Z \|\Lambda(z) \| \d \pi(z) < +\infty, 
\end{equation}where $\| \cdot\| $ is the Euclidean norm. Assume moreover that $1\in \span(\{\Lambda_1,\dots,\Lambda_N\})$. Then, there exist $1\leq K \leq N$, distinct points $z_1,\dots, z_K \in E$ and positive weights $t_1, \dots, t_K$ such that 
\begin{equation}
    \int_Z \Lambda_i(z) \d \pi(z) = \sum^K_{k=1} t_k \Lambda_i(z_k), i=1,\dots,N \quad  \text{ and } \quad \sum^K_{k=1} t_k = \pi(Z).
\end{equation}
\end{thr}
\begin{thr}\label{thm:atomic-minimizer} If the infimum in \eqref{eq:moment-problem} is finite, there exist $1\leq K \leq n+3$, distinct points $z_1,\dots,z_K$ and positive weights $t_1,\dots,t_K$ such that $\pi = \sum^K_{i=1} t_i \delta_{z_i}$ is a minimizer of \eqref{eq:moment-problem}.
\end{thr}

\begin{proof} By \Cref{thm:existence}, there exists a minimizer $\bar \pi \in  \Pi^n_A$ with finite cost. Thus, $\bar \pi$ is concentrated on the set $E= \{z \in Z: c(z) < +\infty, \Theta(z) < +\infty\}$. Define $\tilde c$ and $\tilde \Theta$ to agree with $c$ and $\Theta$ on $E$ and to equal zero on $E^c$. Then define the Borel map $\Lambda : Z \to \R^{n+3}$ by $\Lambda(z) = (w_1(z),\dots, w_n(z), \tilde \Theta(z), \tilde c(z), 1)$. All components $\Lambda_i$ of $\Lambda$ are $\bar \pi$-integrable, so \[ \int \|\Lambda(z) \| \d \bar \pi(z) <+\infty.\]Applying \Cref{thm:tchakaloff} yields a discrete measure that is again a minimizer, since we included $c$ as one of the coordinates of $\Lambda$.
\end{proof}

The set $\Pi^n_A$ is convex, as it is determined by linear constraints. We therefore denote by $\textup{ex}(\Pi^n_A)$ the set of extreme points: if $\pi \in \textup{ex}(\Pi^n_A)$ and $\pi = (1-t)\pi^+ + t(\pi^-)$ for some $\pi^+,\pi^- \in \Pi^n_A$ and $t\in (0,1)$, then $\pi^+=\pi^-$. We also introduce the following notation:
\begin{equation}
    \Phi = (w_1,\dots,w_n,\Theta,1).
\end{equation}
\begin{lemma}\label{lem:moment-perturbation}
Let $\pi\in\textup{ex}(\Pi_A^n)$. If $h\in L^\infty(\pi)$ satisfies
\[
    \int h\Phi\d\pi = 0 
\]
then $h=0$ $\pi$-almost everywhere.
\end{lemma}
\begin{proof} For $0<\epsilon < \frac{1}{1+ \|h\|_\infty}$ define the nonnegative measures $\pi^{\pm} = (1\pm \epsilon h)\pi$. From $\int h \Phi \d\pi = 0$ we have that $\pi^\pm \in \Pi^n_A$, moreover $\pi = \frac{1}{2}(\pi^+ + \pi^-)$. By extremality we have $\pi^+ = \pi^-$ and thus $h \equiv 0$ $\pi$-a.e.
\end{proof}

\begin{thr}\label{thm:extreme-points} \[\textup{ex}(\Pi^n_A) \subset \left\{\begin{aligned} \pi \in \Pi^n_A \, : \, & \pi = \sum^k_{i=1} t_i \delta_{z_i} \text{ for } \{z_i\}^k_{i=1} \subset Z, \sum^k_{i=1} t_i = m, 1\leq k \leq n+2, \text{  } \\&(w_1(z_i),\dots, w_n(z_i),\Theta(z_i),1)^k_{i=1} \text{ are linearly independent}\end{aligned}\right\}\]
\end{thr}
\begin{proof} Let $\pi \in \textup{ex}(\Pi^n_A)$ and assume that $\# \supp(\pi) > n+2$. We can choose $n+3$ pairwise disjoint Borel sets $B_1,\dots B_{n+3}$, each having positive $\pi$-mass, and define the vectors
\begin{equation}
    v_j = \int_{B_j} \Phi \d \pi \in \R^{n+2}, \, \,  j=1\dots, n+3.
\end{equation}As we have $n+3$ vectors $v_j$ in $\R^{n+2}$, these vectors are therefore linearly dependent. Thus, there exists $(a_j)^{n+3}_{j=1}$ not all zero such that
\begin{equation}
    \sum^{n+3}_{j=1} a_j v_j = 0.
\end{equation}Define the function $h = \sum_{j=1}^{n+3} a_j 1_{B_j}$, clearly $h \in L^\infty(\pi)$ and $h \not \equiv 0$. On the other hand, note that
\begin{equation}
    \int h\Phi \d \pi = \sum^{n+3}_{j=1} a_j \int_{B_j}\Phi \d \pi = \sum^{n+3}_{j=1} a_j v_j = 0
\end{equation}and by \Cref{lem:moment-perturbation} we have a contradiction.

Therefore, there exist distinct points $z_1,\dots, z_k \in Z$ and positive weights $t_1,\dots, t_k \in \R$ such that $\pi = \sum^{k}_{i=1} t_i \delta_{z_i}$ for $1\leq k \leq n+2$. Let us show now that the vectors $\Phi(z_1),\dots, \Phi(z_k)$ are linearly independent. Assume that there exist $b_1,\dots,b_k$, not all zero, such that $\sum^k_{i=1} b_i \Phi(z_i) = 0$. Define $h(z_i) = b_i/t_i$ and $h =0$ elsewhere, we have that $h\not \equiv 0$ and $h\in L^\infty(\pi)$. Moreover,
\begin{equation}
    \int h \Phi \d\pi = \sum^k_{i=1} b_i \Phi(z_i) =0
\end{equation}hence, by \Cref{lem:moment-perturbation} we again have a contradiction.
\end{proof}

Since $\C$ is linear and $\Pi^n_A$ is convex and narrowly compact, the minimum of $\C$ over $\Pi^n_A$ is attained at an extreme point, as made precise in the following proposition\corr{ (see \cite{Wink1988MOR, Pine2016MMOR} for the extreme-point approach to optimization over moment sets)}.

\begin{prop}\label{thm:extreme-minimizer} If the infimum in \eqref{eq:moment-problem} is finite, then it is attained at an extreme point of $\Pi^n_A$. In particular, by \Cref{thm:extreme-points}, there exists a minimizer supported on at most $n+2$ points. Moreover, every minimizer is the narrow limit of a sequence of finitely supported minimizers.
\end{prop}

\begin{proof} Define the set 
\begin{equation}F = \{\pi \in \Pi^n_A \, : \, \C(\pi) = \inf_{\pi \in \Pi^n_A} \C(\pi)\}.\end{equation}
$F$ is nonempty by \Cref{thm:existence}. The set $F$ is closed as it is the same as $\{\pi \in  \Pi^n_A \, : \, \C(\pi) \leq \inf_{\pi \in \Pi^n_A} \C(\pi)\}$ and $\C$ is lower semicontinuous; moreover, $F$ is narrowly compact as $\Pi^n_A$ is narrowly compact (\Cref{thm:existence}). 
Moreover, $F$ is convex since $\Pi^n_A$ is convex and $\C$ is linear. By the Krein--Milman theorem, \begin{equation}\label{eq:thr2.6-1}F= \overline{\textup{co}}(\textup{ex}(F))\end{equation}and, as $F\neq \emptyset$, $\textup{ex}(F)\neq \emptyset$ and hence \[\inf_{\pi \in \Pi^n_A }\C(\pi) = \inf_{\pi \in \text{ex}(F)}\C(\pi).\] Now we show that $\text{ex}(F) \subset \text{ex}(\Pi^n_A)$ by showing that $F$ is a face of $\Pi^n_A$. Let $\pi \in F$, $\pi_1,\pi_2\in \Pi^n_A$ such that $\pi=t\pi_1 + (1-t)\pi_2$ with $t\in (0,1)$, then $\C(\pi_1) \geq \inf_{\pi \in \Pi^n_A} \C(\pi)$ and $\C(\pi_2) \geq \inf_{\pi \in \Pi^n_A} \C(\pi)$ and also \[ \C(\pi) =t\C(\pi_1)+(1-t)\C(\pi_2)\] guarantees that $\C(\pi_1) = \C(\pi_2) = \inf_{\pi \in \Pi^n_A} \C(\pi)$, thus $\pi_1,\pi_2 \in F$ which shows that $F$ is a face.

Lastly, let $\pi \in F$. The identity \eqref{eq:thr2.6-1} gives a sequence $\pi_k \in \operatorname{co}(\operatorname{ex}(F))$ such that $\pi_k \weakto \pi$ narrowly, $\pi_k$ is a finite convex combination of extreme points of $F$, which are finitely supported, thus so is $\pi_k$, and finally $\pi_k \in F$, so every $\pi_k$ is a minimizer.
\end{proof}

\subsection{Duality}\label{sec:duality}

In what follows, we consider $\Theta: Z \to [0,+\infty)$ continuous and the $\Theta$-weighted vector space
\[\M_\Theta(Z)
\coloneqq
\left\{
\pi\in \M(Z):\int_Z \Theta\,d|\pi|<+\infty
\right\},\]
where \(\M(Z)\) denotes the space of finite signed Radon measures on \(Z\), and
\(|\pi|\) denotes the total variation measure of \(\pi\). We also set
\[\C_{b,\Theta}(Z)
\coloneqq
C_b(Z)+\operatorname{span}\{\Theta\}.\]
The pairing between \(\M_\Theta(Z)\) and \(\C_{b,\Theta}(Z)\) is given by
\[
\langle \pi,u\rangle
\coloneqq
\int_Z u\,d\pi,
\qquad
\pi\in \M_\Theta(Z),\ u\in \C_{b,\Theta}(Z).
\]

We endow \(\M_\Theta(Z)\) with the weak topology $\sigma(\M_\Theta(Z),\C_{b,\Theta}(Z)),$ i.e.\ the topology  generated by the maps
\[\pi\mapsto \int_Z f\,d\pi,
\, f\in C_b(Z) \text{ and } \pi\mapsto \int_Z \Theta\,d\pi.\]
By
\cite[Theorem 3.10]{rudin1991functional}, we have
$\left(\M_\Theta(Z),\sigma(\M_\Theta(Z),\C_{b,\Theta}(Z))\right)^*
\simeq
\C_{b,\Theta}(Z).$

For $m> 0$, $A > 0 $ and $W = \{w_1,\dots,w_n\}$  consider the functional $ \D: C_b(Z) \times \R \times \R_+ \to \R\cup\{-\infty\}$ \begin{equation}
    \D(w,\sigma,\tau) = \begin{cases}
        \sigma m - \tau A \text{ if } (w,\sigma, \tau) \in \Gamma^n_\Theta\\
        -\infty \text{ otherwise}
    \end{cases}\end{equation}
where \[\Gamma^n_\Theta\coloneq \{(w,\sigma,\tau) \in C_b(Z) \times \R \times \R_+ \, | \, w \in \textup{cone}_+(W), \tau \geq 0, \text{ and }  w+\sigma \leq c+\tau \Theta \}.\]

\begin{thr}\label{thm:duality} Let $Z$, $c$, $\Theta$ and $W= \{w_1,\dots,w_n\}$ be as in \eqref{eq:moment-problem}. Assume that there exists $\eta \in \Pi^n_A$ such that  $\int c \d\eta <+\infty$ and $\supp(\eta)= \overline{\{c<+\infty\}}$. Then the following duality holds \[\max_{(w,\sigma,\tau)} \mathcal \D(w,\sigma,\tau) = \min_{\pi \in \Pi^n_A} \C(\pi).\] Moreover, for any minimizer $\pi^\star$ of \eqref{eq:moment-problem} and dual maximizer $(w^\star,\sigma^\star,\tau^\star)$, the following optimality condition holds:
\begin{equation}
    c(z) + \tau^\star\Theta(z) = w^\star(z) + \sigma^\star \qquad \pi^\star\text{-almost everywhere}.
\end{equation}
\end{thr}

    \begin{proof} Recall that on  $\M_\Theta(Z)$ we consider the $\Theta$-weighted topology. Define $\mathcal F : \mathcal M_\Theta(Z) \to \R \cup\{+\infty\}$ by \[\mathcal F (\pi) = \begin{cases}
            \int_Z c \d \pi \text{ if } \pi \text{ a positive measure} \\
            +\infty  \text{ otherwise}
        \end{cases}\] note that $\mathcal F $ is proper (since it is finite on $\pi=0$) and lower semicontinuous, as $c$ is lower semicontinuous and bounded below, and $\mathcal F$ is convex. 

    Define $\Xi = \{(a,s,t)\in \R^n \times \R \times \R \, | \, a_i \geq 0, \text{ for } 1\leq i \leq n, \, s=m,\ t \leq A\}$ and let $\mathcal G :\R^{n+2} \to \R \cup \{+\infty\}$ be given by \[\mathcal G(a,s,t) = \begin{cases} 0 \text{ if } (a,s,t) \in  \Xi \\ +\infty \text{ otherwise }
    \end{cases}\]note that $\Xi$ is non-empty, closed and convex, thus $\mathcal G$ is proper, lower semicontinuous and convex. Moreover, $\Xi$ is a polyhedron, so $\mathcal G$ is polyhedral.

    Finally, define the linear transformation $L:\mathcal M_\Theta(Z) \to \R^{n+2}$ as \[L(\pi) = \left( \int_{Z} w_1 \d \pi, \dots, \int_{Z} w_n \d \pi, \int_Z1\d\pi, \int_{Z} \Theta \d\pi \right)\] and note that $L$ is continuous for the $\Theta$-weighted topology.

    We have
    \begin{equation}\label{eq:duality-primal}\inf_{\pi \in \mathcal{M}_\Theta(Z)} [\mathcal F(\pi) + \mathcal G(L\pi)] = \inf_{\pi \in \Pi_A^n} \int c\d\pi.\end{equation}


 By assumption, there exists $\eta$ such that $L(\eta) \in \Xi$. It remains to check that $L(\eta) \in \textup{ri}(L(\text{dom}(\mathcal F)))$, where the relative interior is taken with respect to $H\coloneq \text{aff}(L(\textup{dom}(\mathcal F)))$. Since $0\in H$, $H$ is a linear space. Thus, we can form a basis of $H$ by choosing points $z_1,\dots,z_q \in \{c<+\infty\}$ and considering the vectors $L(\delta_{z_1}),\dots, L(\delta_{z_q})$. Choose pairwise disjoint neighborhoods $U_i$ of $z_i$; the vectors $L(\eta|_{U_1}),\dots,L(\eta|_{U_q})$ also form a basis of $H$, since they are nonzero (here we use the full-support assumption on $\eta$) and linear independence is preserved under small perturbations. Consider the measure \[\eta_t = \eta + \sum^q_{j=1} t_j \eta |_{U_j},\]where each $t_j$ varies in a small neighborhood of zero. For instance, taking $|t_j|<1/2$ ensures that $\eta_t \in \text{dom}(\mathcal F)$. By linearity, \[L(\eta_t) = L(\eta) + \sum^q_{j=1}t_jL(\eta |_{U_j}).\]Because the vectors in the sum form a basis, varying each $t_j$ produces a neighborhood of $L(\eta)$ in $H$. Hence, $L(\eta) \in \text{ri}(L(\text{dom}(\mathcal F)))$, as required. Since $\mathcal G$ is polyhedral, the relative-interior condition above allows us to apply Fenchel--Rockafellar duality (see \cite[Theorem 4.2 (ii)]{borwein1992partially}), and we have
    \begin{equation}\label{eq:fenchel-duality}\inf_{\pi \in \M_\Theta(Z)} [\mathcal F(\pi) + \mathcal G(L\pi)] = \max_{x^*} [-\mathcal G^*(x^*) - \mathcal F^*(-L^*x^*)].\end{equation}  
 
    Let us compute $\mathcal F^*: \mathcal M_\Theta(Z)^* \to \R \cup \{+\infty\}$ and $\mathcal G^*: \R^{n+2} \to \R\cup\{+\infty\}$. Let $h \in \mathcal M_\Theta(Z)^*$, by definition of convex conjugate we have \begin{equation}\label{eq:cost-conjugate}\begin{aligned}\mathcal F^*(h) &= \sup_{\pi\in \M_\Theta(Z)\cap \mathcal M_+(Z)} \langle \pi,h\rangle - \mathcal F (\pi) \\ &= \sup_{\pi \in \M_\Theta(Z) \cap \mathcal M_+(Z)}\int_{Z} h-c \d \pi \\ &= \begin{cases}
        0 \text{ if } h\leq c \\ +\infty \text{ otherwise}  
    \end{cases}\end{aligned}\end{equation}

    For $\mathcal G^*$ let $(\alpha, \sigma, \tau) \in \R^{n+2}$, we have
    \begin{equation}\label{eq:constraint-conjugate}\begin{aligned}\mathcal G^*(\alpha,\sigma, \tau ) &= \sup_{(a,s,t) \in \R^{n+2}} \langle (a,s,t),(\alpha,\sigma,\tau)\rangle - \mathcal G(a,s,t) \\ &= \sup_{\substack{a_i \geq 0, \\
    s=m \\
    t\leq A}}\sum^n_{i=1}\alpha_ia_i + \sigma s  +  \tau t \\
    &=\begin{cases}
       m\sigma + A\tau  \text{ if } \alpha_i \leq 0 \text{ for } 1\leq i \leq n,  \tau \geq 0 \\
        + \infty \text{ otherwise}
    \end{cases}.
    \end{aligned}\end{equation} Finally, the adjoint of $L$, $L^*$, is computed through the identity $\langle (L^*(\alpha,\sigma,\tau), \pi)  \rangle = \langle (\alpha, \sigma,\tau), L\pi \rangle$ to be \begin{equation}\label{eq:moment-adjoint}L^*(\alpha,\sigma, \tau) = \sum^n_{i=1} \alpha_iw_i + \sigma  + \tau \Theta .\end{equation}

    Substituting \eqref{eq:constraint-conjugate}, \eqref{eq:cost-conjugate}, and \eqref{eq:moment-adjoint} into \eqref{eq:fenchel-duality} we obtain, by using that $\mathcal F^*(-L^*x^*) = 0$ when $-\sum^n_{i=1} \alpha_iw_i - \sigma  - \tau \Theta \leq c$,
    \begin{equation}\label{eq:moment-dual}
    \begin{aligned}
        \sup_{x^*} -\mathcal G^*(x^*) &- \mathcal F^*(-L^*x^*) =  \\ &=\sup_{\substack{(\alpha,\sigma,\tau)\in \R^{n+2} \\ \tau \geq 0, \alpha_i \leq 0 \\-\sum^n_{i=1} \alpha_iw_i \leq c + \tau \Theta+\sigma  }} - \sigma m - \tau A    \\
        & = \sup_{\substack{(w,\sigma,\tau ) \in \text{cone}_+(W) \times \R \times \R_+ \\ w+\sigma  \leq c+ \tau \Theta  }}  \sigma m -\tau A 
    \end{aligned}\end{equation}where we absorbed the negative signs into the coefficients $\alpha_i$ and $\sigma$ 
    and the conclusion follows from \eqref{eq:duality-primal}.

    Let us now show the optimality condition. Take $\pi^\star \in \Pi_A^n$ to be a minimizer of \eqref{eq:moment-problem} and $(w^\star,\sigma^\star,\tau^\star)$ a maximizer of $\D$. By strong duality we have that 
    \begin{equation}
        \begin{aligned}
            0 = \int c \d \pi^\star - \sigma^\star m + \tau^\star A = \int (c - \sigma^\star) \d \pi^\star  + \tau^\star A
        \end{aligned}
    \end{equation}
    adding and subtracting the terms $\int w^\star \d \pi^\star $ and $\tau^\star\int \Theta \d \pi^\star$ we obtain
    \begin{equation}
        \begin{aligned}
            0 = \int (c + \tau^\star \Theta - w^\star - \sigma^\star) \d \pi^\star  +\int w^\star \d \pi^\star + \tau^\star \left (A - \int \Theta\d\pi^\star\right)
        \end{aligned}
    \end{equation}and, as every term on the right-hand side is nonnegative, it follows that \begin{equation}
        c + \tau^\star \Theta - w^\star - \sigma^\star = 0 \, \pi^\star\text{-a.e.}\text{, } \int w^\star \d \pi^\star = 0 \text{ and } \tau^\star \left (A - \int \Theta\d\pi^\star\right) = 0.
    \end{equation}
    \end{proof}

\subsection{\texorpdfstring{$\Gamma$}{Gamma}-convergence}\label{sec:gamma-convergence}

Thus far we have kept the set of moment functions $\{w_1,\dots,w_n\}$ fixed. In this subsection, we want to let $n\to \infty$ and determine the limit problem of \eqref{eq:moment-problem}. Let us start by fixing some notation. For $n\in \N$ let $\Pi^n_A$ be the admissible set, as defined in \eqref{eq:moment-set} (these sets need not be nested), define $\Pi^\infty_A = \bigcap_{n\in \N} \Pi^n_A$ and the following functionals
$\C^n,\C^\infty : \mathcal M_+(Z)\to \R\cup\{+\infty\}$ by
\[
\C^n(\pi) =
\begin{cases}
\displaystyle \int_{Z} c \,\d\pi & \text{if } \pi \in \Pi_A^n, \\[0.2cm]
+\infty & \text{otherwise,}
\end{cases}
\qquad
\C^\infty(\pi) =
\begin{cases}
\displaystyle \int_{Z} c \,\d\pi & \text{if } \pi \in \Pi_A^\infty, \\[0.2cm]
+\infty & \text{otherwise.}
\end{cases}
\]

\begin{thr}\label{thm:moment-gamma}
Assume that whenever $n_k\to+\infty$ and
    $\pi_k\in\Pi_A^{n_k}$ satisfy $\pi_k\weakto\pi$ narrowly,
    we have $\pi\in\Pi_A^\infty$.
Then $\C^n$ $\Gamma$-converges to $\C^\infty$ with respect
to narrow convergence, and the sequence $(\C^n)_{n\in\N}$
is equicoercive.
\end{thr}

In the Monge--Kantorovich case, the assumption of \Cref{thm:moment-gamma} holds when the families of moment functions are increasing and fully characterize the marginals; see \Cref{thm:ot-convergence}.

\begin{proof}
We start by proving the $\Gamma$-$\liminf$ inequality.
Let $\pi\in\M_+(Z)$ and consider $\pi^n\weakto\pi$ narrowly.
If $\liminf_{n\to\infty}\C^n(\pi^n)=+\infty$, there is nothing
to prove. Assume therefore that
$\liminf_{n\to\infty}\C^n(\pi^n)<+\infty$.
Then we may extract a subsequence $(n_k)_{k\in\N}$ such that
$\pi^{n_k}\in\Pi_A^{n_k}$ for every $k\in\N$ and
\[
    \liminf_{n\to\infty}\C^n(\pi^n)
    =
    \lim_{k\to\infty}\C^{n_k}(\pi^{n_k}).
\]

Let us show that this implies $\pi\in\Pi_A^\infty$.
Since $\pi^{n_k}\in\Pi_A^{n_k}$ and
$\pi^{n_k}\weakto\pi$ narrowly, this follows directly from
the assumption.

Finally, by lower semicontinuity of $c$ we obtain
\[
    \C^\infty(\pi)
    =\int_Z c\,\d\pi
    \leq\liminf_{k\to\infty}\int_Z c\,\d\pi^{n_k}
    =\liminf_{n\to\infty}\C^n(\pi^n).
\]

For the $\Gamma$-$\limsup$ inequality, consider now $\pi\in\M_+(Z)$.
If $\pi\notin\Pi_A^\infty$, then the limsup inequality is
trivially true. Assume $\pi\in\Pi_A^\infty$ and consider
the sequence $(\pi^n)_{n\in\N}$ defined by $\pi^n=\pi$
for every $n\in\N$.

By definition of $\Pi^\infty_A$, $\pi\in\Pi_A^n$ for every $n$,
so $\C^n(\pi)=\C^\infty(\pi)$ for every $n$. Therefore
\[
    \limsup_{n\to\infty}\C^n(\pi^n)
    =\limsup_{n\to\infty}\C^n(\pi)
    =\C^\infty(\pi).
\]

Finally, let us show that $\{\C^n\}_{n\in\N}$ is equicoercive.
By Prokhorov's theorem and the fixed mass constraint, it suffices
to show that, for every $M>0$, the sublevel sets
$\{\pi:\C^n(\pi)\leq M\}$ are uniformly tight with respect to $n$.
Fix $R>0$. If $\C^n(\pi)\leq M<+\infty$, then
$\pi\in\Pi_A^n$. Note that, since
$\Theta:Z\to[0,+\infty]$ is coercive, the set
$K\coloneq\{z:\Theta(z)\leq R\}$ is compact. Hence,
\[
\begin{aligned}
    \pi(\{\Theta>R\})
    &=\int_Z\chi_{\{\Theta>R\}}\,\d\pi(z) \leq\frac{1}{R}\int_Z\Theta(z)\,\d\pi(z) \leq\frac{A}{R},
\end{aligned}
\]
where the last inequality follows from $\pi\in\Pi_A^n$.
The result follows by observing that the right-hand side
converges to $0$ as $R\to+\infty$, uniformly over all these
sublevel sets.
\end{proof}

\subsection{The Monge--Kantorovich case}\label{sec:mk}

Let $X$ and $Y$ be Polish spaces, let $\mu\in\P(X)$ and
$\nu\in\P(Y)$, and let
$c:X\times Y\to\R\cup\{+\infty\}$ be lower semicontinuous and
bounded below. The Monge--Kantorovich problem is
\begin{equation}\label{eq:ot}
    \inf_{\pi\in\Pi(\mu,\nu)}
    \int_{X\times Y} c(x,y)\,\d\pi(x,y),
\end{equation}
where
\[
    \Pi(\mu,\nu)
    =
    \bigl\{
    \pi\in\P(X\times Y):
    (p_X)_\#\pi=\mu,\ (p_Y)_\#\pi=\nu
    \bigr\},
\]
and $p_X:X\times Y\to X$ and $p_Y:X\times Y\to Y$ denote the
canonical projections.

Choose nonnegative, lower semicontinuous functions $ 
    \theta_X:X\to[0,+\infty],$ $ 
    \theta_Y:Y\to[0,+\infty]$
with compact sublevel sets and finite integrals against $\mu$
and $\nu$, respectively. Set
\[
    \Theta(x,y)=\theta_X(x)+\theta_Y(y)
\]
and choose $A>0$ such that
\[
    A\geq
    \int_X\theta_X\,\d\mu+\int_Y\theta_Y\,\d\nu.
\]

We set $Z = X\times Y$ and $m=1$. As moment functions, we choose $\varphi_0,\varphi_1,\dots,\varphi_r \in C_b(X)$ and $\psi_0,\psi_1,\dots,\psi_s \in C_b(Y)$ where $\varphi_0 \equiv \psi_0 \equiv 1$. Let the functions $w_1,\dots,w_n$ of \eqref{eq:moment-problem} be the functions $f_i,-f_i,g_j,-g_j$, $0\leq i\leq r$, $0\leq j\leq s$, where
\[
    f_i=\varphi_i\circ p_X-\int_X\varphi_i\,\d\mu,
    \qquad
    g_j=\psi_j\circ p_Y-\int_Y\psi_j\,\d\nu,
\]
and we take $\Theta$ and $A$ as above. We define $\Pi^n_A(\mu,\nu)$ to be $\Pi^n_A$ with the choices above and we consider the problem
\begin{equation}\label{eq:moment-ot}
    \inf_{\pi \in \Pi^n_A(\mu,\nu)} \int_{X\times Y} c(x,y) \d \pi(x,y).
\end{equation}

\begin{thr}
\label{thm:ot-discrete}
If \eqref{eq:moment-ot} is finite, then it admits a minimizer which is
an extreme point of $\Pi_A^n(\mu,\nu)$ and it is supported on at most $r+s+2$ points.
\end{thr}

This follows from \Cref{thm:extreme-points} and \Cref{thm:extreme-minimizer}: up to signs and additive constants, the components of $\Phi$ are $\varphi_1\circ p_X,\dots,\varphi_r\circ p_X$, $\psi_1\circ p_Y,\dots,\psi_s\circ p_Y$, $\Theta$ and $1$, so that at most $r+s+2$ of the vectors $\Phi(z_i)$ can be linearly independent.

For the duality result below, we introduce the following notation: $E_r = \span(\{\varphi_0,\dots,\varphi_r\})$ and $E_s = \span(\{\psi_0,\dots,\psi_s\})$.

\begin{thr}
\label{thm:ot-duality}
Assume that $\Theta$ is finite-valued and continuous.
If there exists $\eta\in\Pi_A^n(\mu,\nu)$ such that
$
    \int_{X\times Y}c\,\d\eta<+\infty$ and $
    \supp(\eta)=\overline{\{c<+\infty\}},
$
then
\begin{equation}\label{eq:ot-dual}
\begin{aligned}
    \min_{\pi\in\Pi_A^n(\mu,\nu)}
    \int_{X\times Y}c\,\d\pi
    =
    \max_{\substack{
        \varphi\in E_r,\ \psi\in E_s,\ \tau\geq0\\
        \varphi(x)+\psi(y)\leq c(x,y)+\tau\Theta(x,y)
    }}
    \left\{
        \int_X\varphi\,\d\mu
        +\int_Y\psi\,\d\nu-\tau A
    \right\}.
\end{aligned}
\end{equation}
The pointwise inequality is required for every $(x,y)\in X\times Y$.
\end{thr}

\begin{proof} Since $\operatorname{cone}_+(\{w_1,\dots,w_n\}) = \span(\{f_0,\dots,f_r,g_0,\dots,g_s\})$ we can write any $w\in \operatorname{cone}_+(\{w_1,\dots,w_n\})$ as 
\[
\begin{aligned}
    w&=\tilde \varphi\circ p_X-\int_X\tilde \varphi\,\d\mu
       +\tilde \psi\circ p_Y-\int_Y\tilde \psi\,\d\nu,
\end{aligned}
\]for some $\tilde \varphi \in E_r$ and $\tilde \psi \in E_s$. Define
\[\varphi = \tilde \varphi - \int \tilde \varphi \d\mu + \frac{\sigma}{2} \qquad \text{and} \qquad \psi = \tilde \psi - \int \tilde \psi \d \nu + \frac{\sigma}{2} \] to get \[\sigma = \int \varphi \d\mu + \int \psi \d\nu\]and
the conclusion follows from \Cref{thm:duality}. 
\end{proof}
\begin{rmk}[A finite dimensional problem]\label{rmk:finite-dim-dual}
\corr{Writing $\varphi=\sum_{i=0}^r a_i\varphi_i$ and $\psi=\sum_{j=0}^s b_j\psi_j$, the dual problem in \eqref{eq:ot-dual} is a finite-dimensional maximization over the coefficients $(a,b,\tau)\in\R^{r+1}\times\R^{s+1}\times[0,+\infty)$, namely
\[
\begin{aligned}
    \max_{(a,b,\tau)}\ &\sum_{i=0}^r a_i\hat\mu_i+\sum_{j=0}^s b_j\hat\nu_j-\tau A\\
    \text{subject to}\ &\sum_{i=0}^r a_i\varphi_i(x)+\sum_{j=0}^s b_j\psi_j(y)\leq c(x,y)+\tau\Theta(x,y)
    \quad\text{for every }(x,y)\in X\times Y,
\end{aligned}
\]
where $\hat\mu_i=\int_X\varphi_i\,\d\mu$ and $\hat\nu_j=\int_Y\psi_j\,\d\nu$. Hence the marginals enter the dual problem only through their moments, and \eqref{eq:ot-dual} is a semi-infinite linear program, with finitely many unknowns and infinitely many constraints.}
\end{rmk}
\corr{The preceding results concern a fixed finite family of moment functions. We now study convergence to the original transport problem as this family is enlarged. To recover the prescribed marginals in the limit, we require that the full collection of moments determines them uniquely. Following \cite{AlfCoyEhrLom2021Approx}, we express this requirement through the following notion of full characterization.}

\begin{df}\label{def:moment-determinacy} Let $\mu \in \M_+(X)$. A family of  functions $\{\phi^i\}_{i\in \N}$ in $L^1_\mu(X)$ fully characterizes the measure $\mu$ if, for any measure $\eta \in \mathcal{M}_+(X)$ satisfying\[\int_X \phi^i \d \eta = \int_X \phi^i \d \mu, \text{ for all } i \in \N,\] it follows that $\eta = \mu$.
\end{df}

\begin{thr}\label{thm:ot-convergence} Let $(\varphi_i)_{i\in \N} \subset C_b(X)$ and $(\psi_j)_{j\in \N}\subset C_b(Y)$ fully characterize $\mu$ and $\nu$, respectively. Let $I^n$ be the minimum in \eqref{eq:moment-ot} with moment functions $(\varphi_i)^n_{i=1}$ and $(\psi_j)^n_{j=1}$, and let $I$ be the minimum in \eqref{eq:ot}. If $I$ is finite, then $I^n \nearrow I$ as $n\to +\infty$. Moreover, if $\pi^n$ minimizes \eqref{eq:moment-ot} for every $n\in \N$, then, up to a subsequence, $\pi^n$ converges narrowly to a minimizer $\pi$ of \eqref{eq:ot}.
\end{thr}
\begin{proof}
By \Cref{thm:moment-gamma} and the fundamental theorem of $\Gamma$-convergence, it suffices to show that $\Pi^\infty_A = \Pi(\mu,\nu)$ and that the assumption of \Cref{thm:moment-gamma} holds. The families of moment functions being increasing, the sets $\Pi^n_A(\mu,\nu)$ are nested, which also gives the monotonicity of $(I^n)_{n\in\N}$. Moreover, these sets are narrowly closed: if $\pi_k\in\Pi^{n_k}_A(\mu,\nu)$ with $n_k\to+\infty$ converge narrowly to $\pi$, then $\pi\in\Pi^N_A(\mu,\nu)$ for every $N$, that is, $\pi\in\Pi^\infty_A$. The choice of $\Theta$ and $A$ guarantees that
$\Pi(\mu,\nu)\subseteq\Pi_A^n(\mu,\nu)$ for every $n\in\N$.
Conversely, let $\pi\in\bigcap_{n\in\N}\Pi_A^n(\mu,\nu)$.
Then, for every $i,j\in\N$,
\[
    \int_{X\times Y}\varphi_i(x)\,\d\pi
    =\int_X\varphi_i\,\d\mu,
    \qquad
    \int_{X\times Y}\psi_j(y)\,\d\pi
    =\int_Y\psi_j\,\d\nu.
\]
Since these families fully characterize $\mu$ and $\nu$,
we obtain $(p_X)_\#\pi=\mu$ and $(p_Y)_\#\pi=\nu$. Therefore,
\[
    \Pi_A^\infty
    =\bigcap_{n\in\N}\Pi_A^n(\mu,\nu)
    =\Pi(\mu,\nu).
\]
\end{proof}

\corr{Increasing the number of moment constraints does not by itself ensure that the prescribed marginals are recovered in the limit. As the following example shows, two distinct probability measures may have the same moments against every function in an infinite family.}

\begin{exe}\label{ex:indeterminate-moments} Consider $X=[-1,1]$ and the moment functions $\phi_0(x) = 1$ and $\phi_k(x) = x^{2k-1}$ for $k\in \N$. We have that
\[\int \phi_k \d (\delta_0) = \int \phi_k \d \big(\frac{1}{2} (\delta_{-1} + \delta_{1})\big) \text{ for every } k\in \N \]and the measures are not equal.
\end{exe}

\corr{The following lemma provides a sufficient condition for full characterization.}

\begin{lemma}\label{lem:determinacy-c0}
Let $X$ be a locally compact Polish space and let
$(\varphi_i)_{i\in\N}\subset C_b(X)$. If
\[
    C_0(X)\subseteq
    \overline{\operatorname{span}\{\varphi_i:i\in\N\}}^{\|\cdot\|_\infty},
\]
then this family fully characterizes every finite nonnegative
Borel measure on $X$.
\end{lemma}

\begin{proof}
Let $\mu,\eta\in\M_+(X)$ satisfy
$\int_X\varphi_i\,\d\mu=\int_X\varphi_i\,\d\eta$ for every
$i\in\N$. Given $f\in C_0(X)$, choose
$g_k\in\operatorname{span}\{\varphi_i:i\in\N\}$ such that
$\|f-g_k\|_\infty\to0$. Since the integrals of $g_k$ against
$\mu$ and $\eta$ coincide, we have
\[
    \left|\int_X f\,\d\mu-\int_X f\,\d\eta\right|
    \leq \bigl(\mu(X)+\eta(X)\bigr)\|f-g_k\|_\infty
    \longrightarrow0.
\]
Thus $\mu$ and $\eta$ agree on every function in $C_0(X)$ and thus
$\mu=\eta$.
\end{proof}

Let us conclude this discussion of \Cref{def:moment-determinacy} by noting that, even when this condition fails, one still obtains an interesting problem, which we call \textbf{optimal transport between equivalence classes}. 
\begin{df}[Equivalence class]
    For a probability measure $\mu\in \P(X)$, we define its equivalence class with respect to a family of moment functions $\I=\span(\{\varphi_i\}_i)$ as
\[[\mu]_\I\coloneq \left\{\eta \in \mathcal{M}_+(X) \, \bigg| \, \int_X \varphi \d \eta = \int_X \varphi \d \mu, \text{ for all } \varphi \in \I \right\}.\] 
\end{df}
\begin{df}[Optimal transport between equivalence classes]
  The  optimal transport between equivalence classes with confinement bound $A$ is the following minimization problem
\begin{equation}\label{eq:class-ot}
\inf_{\substack{\pi \in \Pi([\mu]_\I,[\nu]_\J) \\ \int_{X\times Y}\Theta\,\d\pi\leq A}} \int_{X\times Y} c(x,y) \d\pi(x,y),
\end{equation}  
where $\Pi([\mu]_\I,[\nu]_\J)$ is the set of all couplings between $[\mu]_\I$ and $[\nu]_\J$ for $\I = \span(\{\varphi_i\}_i)$ and $\J = \span(\{\psi_j\}_j)$. In particular, we say that a measure $\pi \in \P(X\times Y)$ is a coupling between $[\mu]_\I$ and 
$[\nu]_\J$ if there exist representatives $\tilde \mu \in [\mu]_\I$ and $\tilde \nu \in [\nu]_\J$ such that $\pi \in \Pi(\tilde \mu,\tilde \nu)$.
\end{df}

Notice that  \eqref{eq:moment-ot} is a special case of \eqref{eq:class-ot} when $\I$ and $\J$ are finite-dimensional. Moreover, \eqref{eq:class-ot} is the limit problem for \eqref{eq:moment-ot} when the moment functions do not fully characterize the marginals.

\begin{cor}\label{cor:ot-duality}
Under the assumptions of \Cref{thm:ot-convergence},
assume that \eqref{eq:ot} is finite and that
\Cref{thm:ot-duality} applies for every $n\in\N$. Then
\[
    \min_{\pi\in\Pi(\mu,\nu)}\int_{X\times Y}c\,\d\pi
    =
    \sup_{\substack{
        \varphi\in C_b(X),\ \psi\in C_b(Y)\\
        \varphi(x)+\psi(y)\leq c(x,y)
    }}
    \left\{
        \int_X\varphi\,\d\mu+\int_Y\psi\,\d\nu
    \right\}.
\]
\end{cor}

\begin{proof}
Denote the dual supremum by $D$. Weak duality gives $D\leq I$.
Let $(\varphi_n,\psi_n,\tau_n)$ be a maximizer in
\eqref{eq:ot-dual} and define
\[
    f_n=\varphi_n-\tau_n\theta_X,
    \qquad
    g_n=\psi_n-\tau_n\theta_Y.
\]
Then $f_n(x)+g_n(y)\leq c(x,y)$ and
\[
    \int_X f_n\,\d\mu+\int_Y g_n\,\d\nu
    =
    I^n+\tau_n\left(
        A-\int_X\theta_X\,\d\mu-\int_Y\theta_Y\,\d\nu
    \right)
    \geq I^n.
\]
The functions $f_n$ and $g_n$ are continuous and bounded
above. Choose $B_n\geq0$ bounding both functions above,
a lower bound $b$ for $c$, and $M_n\geq0$ such that
$B_n-M_n\leq b$. The bounded continuous functions
\[
    \widehat f_n=\max\{f_n,-M_n\},
    \qquad
    \widehat g_n=\max\{g_n,-M_n\}
\]
remain dual admissible, and,
since truncation from below can only increase the integrals,
$D\geq I^n$. Letting $n\to\infty$ gives $D\geq I$,
which concludes the proof.
\end{proof}

\section{Estimates for the \texorpdfstring{$W_\infty$}{W-infty} distance}\label{sec:wp-winf}
As Example \ref{exe:polynomial-geometry} shows, matching finitely many moments can preserve the optimal value while allowing very different marginal supports. We therefore study which geometric information is retained by localized moment constraints, defined below. Our first result is \Cref{lem:moment-localization}, which shows that, under this choice of moment functions, the support of the moment-relaxed marginals stays quantitatively close to the support of the original marginal. We also establish a \(W_p\)--\(W_\infty\) estimate, \Cref{thm:wp-winf-bound}, in the same spirit as \cite{JylhaRajala2016, BoucJimeRaje2007PAMS} with an additive error \(4r_n\), showing how average transport accuracy controls the maximal displacement needed to couple the relaxed marginals.

\begin{assumption}[Localized moment functions]\label{ass:local-moments}
The measure $\mu\in\P(X)$ has compact and connected support and $\nu\in\P(\supp(\mu))$. For fixed $n\in\N$, the nodes $x_1,\dots,x_n\in\supp(\mu)$ and $r_n>0$ are such that the balls $B(x_i,r_n)$ cover $K=\supp(\mu)$, and the moment functions $\varphi_1,\dots,\varphi_n\in C_b(X)$ are nonnegative and form a partition of unity on $K$, subordinate to this cover, that is,
\begin{equation}\label{eq:local-moments}
    \sum_{i=1}^n\varphi_i=1
    \quad\text{on }K,
    \qquad
    \sum_{i=1}^n\varphi_i\leq1
    \quad\text{on }X,
    \qquad
    \supp(\varphi_i)\subset B(x_i,r_n),
    \quad i=1,\dots,n.
\end{equation}
\end{assumption}
Throughout this section, we work under \Cref{ass:local-moments}. The parameter $r_n$ gives the localization scale of the moment functions. For $\lambda\in\P(\supp(\mu))$, we denote by $[\lambda]_n$ the set of probability measures on $X$ having the same integrals as $\lambda$ against $\varphi_1,\dots,\varphi_n$, and we set
\[
    m_\mu(s)\coloneq\inf_{z\in\supp(\mu)}\mu(B(z,s)),\qquad s>0.
\]

We also write $A^r=\{x\in X :d(x,A)\leq r\}$ for the closed
$r$-neighborhood of a set $A$.

\begin{lemma}\label{lem:moment-localization}
Under \Cref{ass:local-moments}, for every $\lambda\in\P(\supp(\mu))$ and $\lambda_n\in[\lambda]_n$,
\[
    W_\infty(\lambda,\lambda_n)\leq2r_n,
    \qquad
    \supp(\lambda_n)\subset \supp(\mu)^{r_n}.
\]
\end{lemma}

\begin{proof}
Set $S=\sum_{i=1}^n\varphi_i$. Since $S=1$ on $\supp(\mu)$,
the moment constraints give
$\int_{X}(1-S)\,\d\lambda_n=0$.
The function $1-S$ is continuous and nonnegative, so it
vanishes on $\supp(\lambda_n)$. Thus, for every
$y\in\supp(\lambda_n)$, there exists $i$ with
$\varphi_i(y)>0$. The support condition on $\varphi_i$
then gives $d(y,x_i)<r_n$, proving the support inclusion.

To estimate the displacement, set
$a_i=\int_{X}\varphi_i\,\d\lambda$ and define
\[
    \gamma(\d x,\d y)
    =
    \sum_{i:a_i>0}
    \frac{\varphi_i(x)\varphi_i(y)}{a_i}
    \lambda(\d x)\lambda_n(\d y).
\]
The moment identities and the equality $S=1$ almost
everywhere for both measures show that
$\gamma\in\Pi(\lambda,\lambda_n)$.
Each summand is supported on
$\supp(\varphi_i)\times\supp(\varphi_i)$, where
$d(x,y)\leq2r_n$. Hence
$W_\infty(\lambda,\lambda_n)\leq2r_n$.
\end{proof}

\begin{lemma}\label{lem:layer-mass}
Assume \Cref{ass:local-moments}. Let $A\subset X$ be closed and let $R>0$.
Suppose that there exist $z_0,z_1\in \supp(\mu)$ such that
\[
    d(z_0,A)\leq r_n,
    \qquad
    d(z_1,A)>R-r_n.
\]
Then, for every $s>0$ satisfying $R>4r_n+3s$ and every
$\mu_n\in[\mu]_n$,
\[
    \mu_n(A^{R-s}\setminus A)\geq m_\mu(s).
\]
\end{lemma}

\begin{proof}
Choose $0<\eta<R-4r_n-3s$.
The number $2r_n+s+\eta$ lies between $r_n$ and $R-r_n$.
Since $\supp(\mu)$ is connected and the distance to $A$ is continuous,
there exists $z\in \supp(\mu)$ such that
$d(z,A)=2r_n+s+\eta$.

By the construction in \Cref{lem:moment-localization}, there
exists $\gamma\in\Pi(\mu,\mu_n)$ with displacement at most
$2r_n$. For $\gamma$-almost every pair $(x,y)$ with
$x\in B(z,s)$, we have
\[\begin{aligned}
    d(y,z)
    &\leq d(y,x)+d(x,z)<2r_n+s,
\\
    |d(y,A)-d(z,A)|
    &\leq d(y,z),
\end{aligned}\]
where the second inequality follows because the distance
to $A$ is $1$-Lipschitz. Consequently,
$\eta<d(y,A)<4r_n+2s+\eta<R-s$, so
$y\in A^{R-s}\setminus A$.

Thus all the mass of $\mu$ in $B(z,s)$ is transported
into $A^{R-s}\setminus A$, and
\[
\begin{aligned}
    \mu_n(A^{R-s}\setminus A)
    &\geq
    \gamma\bigl(B(z,s)\times(A^{R-s}\setminus A)\bigr)\\
    &=\mu(B(z,s))
    \geq m_\mu(s).
\\ \mathstrut
\end{aligned}\qedhere\]
\end{proof}
\vspace{-0.5cm}

\begin{thr}\label{thm:wp-winf-bound}
Under \Cref{ass:local-moments}, for every $1\leq p<\infty$, $s>0$ and $(\mu_n,\nu_n)\in[\mu]_n\times[\nu]_n$,
\[
    W_p^p(\mu_n,\nu_n)<m_\mu(s)s^p
    \quad\Longrightarrow\quad
    W_\infty(\mu_n,\nu_n)\leq4r_n+3s.
\]
\end{thr}

\begin{proof}
Fix $s>0$ such that
$W_p^p(\mu_n,\nu_n)<m_\mu(s)s^p$, and let
$\pi^n\in\Pi(\mu_n,\nu_n)$ be optimal for $W_p$.
The transport cost bounds the mass displaced by at least $s$:
\begin{equation}\label{eq:large-displacement}
\begin{aligned}
    \pi^n(\{d(x,y)\geq s\})
    &\leq
    \frac{1}{s^p}\int_{X\times X}d(x,y)^p\,\d\pi^n\\
    &=\frac{W_p^p(\mu_n,\nu_n)}{s^p}
    <m_\mu(s).
\end{aligned}
\end{equation}

We argue by contradiction. Suppose that
$W_\infty(\mu_n,\nu_n)>4r_n+3s$, and choose
$4r_n+3s<R<W_\infty(\mu_n,\nu_n)$.
By Strassen's theorem (see \cite[Theorem~1.27 and Corollary~1.28]{Villani2003Topics}), there exists a closed set
$A\subset X$ such that
\begin{equation}\label{eq:strassen-separation}
    \mu_n(A)>\nu_n(A^R).
\end{equation}
In particular, $\mu_n(A)>0$ and
$\nu_n((A^R)^c)>0$.
By \Cref{lem:moment-localization}, both measures are supported
in $\supp(\mu)^{r_n}$. We can therefore find $z_0,z_1\in \supp(\mu)$ with
$d(z_0,A)\leq r_n$ and $d(z_1,A)>R-r_n$.
Applying \Cref{lem:layer-mass} gives
\begin{equation}\label{eq:layer-mass}
    \mu_n(A^{R-s}\setminus A)\geq m_\mu(s).
\end{equation}

On the other hand, if $x\in A^{R-s}$ and $y\notin A^R$,
then $d(x,y)>s$. Using the marginals of $\pi^n$ and
\eqref{eq:layer-mass}, we obtain
\[\begin{aligned}
    \mu_n(A^{R-s})
    &\leq
    \nu_n(A^R)+\pi^n(\{d(x,y)\geq s\}),
\\
    \mu_n(A^{R-s})
    &=\mu_n(A)+\mu_n(A^{R-s}\setminus A)
    \geq\mu_n(A)+m_\mu(s).
\end{aligned}\]
Combining these inequalities with \eqref{eq:large-displacement}
yields $\mu_n(A)<\nu_n(A^R)$, contradicting
\eqref{eq:strassen-separation}.
\end{proof}

The next corollary is similar to \cite[Theorem 1.2]{BoucJimeRaje2007PAMS} and to the Euclidean case of \cite[Theorem 1.1]{JylhaRajala2016}. Under the assumptions below, $m_\mu(t)\geq ct^d$ for sufficiently small $t>0$. Let us highlight that both $\mu_n$ and $\nu_n$ are allowed to be discrete measures, while the lower mass bound is imposed on the reference measure $\mu$. Thus, this corollary provides an analogue of these $L^\infty$ estimates for discrete measures satisfying localized moment constraints, with an additive error $4r_n$.

\begin{cor}\label{cor:wp-winf-density} Let $X= \R^d$. 
Under \Cref{ass:local-moments},
 let $\Omega\subset\R^d$  be a bounded, connected
domain with Lipschitz boundary and
$\mu=f\mathbf1_\Omega\,\d x$, where $f\geq a>0$
almost everywhere on $\Omega$.
Then there exists $c_{p,d,a,\Omega}>0$ such that
\[
    W_\infty(\mu_n,\nu_n)
    \leq
    4r_n+c_{p,d,a,\Omega}
    W_p^{\frac{p}{p+d}}(\mu_n,\nu_n)
\]
for every $(\mu_n,\nu_n)\in[\mu]_n\times[\nu]_n$.
\end{cor}

\begin{proof}
Here $\supp(\mu)=\overline\Omega$. Set $D=\operatorname{diam}(\supp(\mu))$.
The Lipschitz boundary condition gives a constant
$c_\Omega>0$ such that
$|\Omega\cap B(z,s)|\geq c_\Omega s^d$
for every $z\in \supp(\mu)$ and $0<s\leq D$.
Since $f\geq a$, it follows that
$m_\mu(s)\geq a c_\Omega s^d$ on this range.

Set
\[
    t=(a c_\Omega)^{-\frac1{p+d}}
    W_p^{\frac p{p+d}}(\mu_n,\nu_n).
\]
If $t<D$, then every $s\in(t,D]$ satisfies
$W_p^p(\mu_n,\nu_n)<a c_\Omega s^{p+d}
\leq m_\mu(s)s^p$.
By \Cref{thm:wp-winf-bound},
$W_\infty(\mu_n,\nu_n)\leq4r_n+3s$.
Letting $s\downarrow t$ gives the desired estimate.

If $t\geq D$, both measures are supported in $\supp(\mu)^{r_n}$,
whose diameter is at most $D+2r_n$. Hence
$W_\infty(\mu_n,\nu_n)\leq D+2r_n\leq4r_n+3t$.
The result follows in both cases with
$c_{p,d,a,\Omega}=3(a c_\Omega)^{-1/(p+d)}$.
\end{proof}

\begin{rmk} In particular, if $r_n \to 0$ as $n\to +\infty$, since $c_{p,d,a,\Omega}$ does not depend on $n$, one recovers the inequality originally proven in \cite{BoucJimeRaje2007PAMS, JylhaRajala2016}. Note that our strategy is similar to Lemma 2.2 in \cite{JylhaRajala2016}, but we obtain an improvement on the underlying constants. 
\end{rmk}

 \section{Estimates for repulsive costs}\label{sec:repulsive-costs}

In \Cref{sec:wp-winf}, the moment functions were localized and the estimates concerned the geometry of the relaxed marginals. We now turn to multi-marginal problems with repulsive costs, such as the Coulomb cost, for which the relevant question is whether minimizers of the relaxed problem stay away from the diagonal, where the cost is singular. Here the moment functions need not be localized. We prove a diagonal exclusion estimate which is uniform in $n$ (\Cref{thm:diagonal-estimate} and \Cref{cor:diagonal-uniform}), and we derive two consequences: the cost can be truncated without changing the minimizers (\Cref{thm:repulsive-truncation}) and, for localized moment functions, the dual problem admits potentials with explicit bounds (\Cref{thm:repulsive-potentials}).

Fix $M\geq2$ and let $X$ be a Polish space endowed with a distance $d$. We consider costs of the form
\begin{equation}\label{eq:repulsive-cost}
    c(x_1,\dots,x_M)=\sum_{1\leq i < j \leq M} h(d(x_i,x_{j})),
\end{equation}
where $h:(0,+\infty)\to[0,+\infty)$ is strictly decreasing and continuous, with $\lim_{r\to0^+}h(r)=+\infty$, and we set $h(0)=+\infty$. Let $\rho\in\P(X)$. For each $n\in\N$, we take the same family $\I^n =\{\varphi_i\}_{i=1}^n\subset C_b(X)$ of moment functions for all $M$ marginals; these moment functions need not form a partition of unity. We denote by $\Pi(\rho)$ the set of probability measures on $X^M$ whose marginals are all equal to $\rho$, and by $\Pi^n(\rho)$ the set of probability measures on $X^M$ whose marginals have the same integrals as $\rho$ against every function in $\I^n$. To guarantee the existence of minimizers, we also impose a confinement constraint.
 Let $\theta:X\to[0,+\infty)$ be continuous with compact sublevel sets, assume $\int_X\theta\,\d\rho<+\infty$, and choose $A>0$ such that $A\geq M\int_X\theta\,\d\rho$. Set $\Theta(x_1,\dots,x_M)=\theta(x_1)+\dots +\theta(x_M)$ and define for every $n\in\N$
\begin{equation}\label{eq:repulsive-moment-set}
    \Pi_A^n(\rho)\coloneq
    \Bigl\{\pi\in\Pi^n(\rho):
    \int_{X^M}\Theta\,\d\pi\leq A\Bigr\},
\end{equation}
\begin{equation}\label{eq:repulsive-moment-problem}
    I^n(\rho)\coloneq
    \inf_{\pi\in\Pi_A^n(\rho)}\C(\pi),
\end{equation}
and
\begin{equation}\label{eq:repulsive-problem}
    I(\rho)\coloneq
    \inf_{\pi\in\Pi(\rho)}\C(\pi),
\end{equation}
where, as before, $\C(\pi)=\int_{X^M}c\,\d\pi$. Problem \eqref{eq:repulsive-moment-problem} is \eqref{eq:moment-problem} with $Z=X^M$, $m=1$ and moment functions $\pm\bigl(\varphi_i\circ p_j-\int_X\varphi_i\,\d\rho\bigr)$, $1\leq i\leq n$, $1\leq j\leq M$, where $p_j$ denotes the projection onto the $j$-th factor.
For $\alpha>0$, we set
\[
    D_\alpha\coloneq
    \bigl\{(x_1,\dots,x_M)\in X^M:d(x_i,x_j)<\alpha, \text{ for some } i\neq j\bigr\},
    \qquad
    \kappa_\rho(\alpha)\coloneq\sup_{x\in X}\rho(B(x,\alpha)).
\]

\begin{prop}\label{prop:finite-cost-separation}
Let $\pi\in\Pi_A^n(\rho)$ satisfy $\int_{X^M}c\,\d\pi<+\infty$. Then there exists a finitely supported sequence $\pi_k \in \Pi^n_A(\rho)$ converging narrowly to $\pi$ such that, for every $k\in \N$,
\[\supp(\pi_k) \subset \supp(\pi), \, \, \#\supp(\pi_k) \leq Mn+k+3,\]
\[\int_{X^M} c \d \pi_k = \int_{X^M} c \d \pi, \, \, \int_{X^M}\Theta \d \pi_k = \int_{X^M} \Theta \d \pi\]and
$\pi_k(D_{\delta_k}) = 0$, where $\delta_k>0$ and $\delta_k \to 0$ as $k\to +\infty$.
\end{prop}

\begin{proof}
Choose a countable family $\{w_j\}_{j\in\N}\subset C_b(X^M)$ fully characterizing $\pi$ in the sense of \Cref{def:moment-determinacy}. The finiteness of the cost implies that $\pi$ is concentrated on the Borel set
$
    E=\supp(\pi)\cap\{c<+\infty\}.
$
Define $\widetilde c$ to agree with $c$ on $\{c<+\infty\}$ and to equal zero elsewhere. For each $k\in\N$, consider the following set of moment functions
\begin{equation}\label{eq:Prop4.1-fcts}
    \bigcup^M_{j=1}\left(\{\varphi_i\circ p_j\}_{i=1}^n\right)
    \cup\{1,\Theta,\widetilde c,w_1,\dots,w_k\}.
\end{equation}
By \Cref{thm:tchakaloff}, there exist $1\leq L_k\leq Mn+k+3$, distinct points $x_{k,\ell}\in E$ and positive weights $a_{k,\ell}$ such that $\pi_k=\sum_{\ell=1}^{L_k}a_{k,\ell}\delta_{x_{k,\ell}}$ and $\pi_k$ has the same integral as $\pi$ with respect to all functions in \eqref{eq:Prop4.1-fcts}, in particular $\pi_k\in\Pi_A^n(\rho)$. Since the points $x_{k,\ell} = (x^1_{k,\ell},\dots, x^M_{k,\ell}) \in E$ it follows that $x^i_{k,\ell}\neq x^j_{k,\ell}$ for $i\neq j$ and every $1\leq \ell\leq L_k$.

As the set $\Pi^n_A$ is narrowly compact, the sequence $(\pi_k)_{k\in \N}$ is tight and, by  full characterization, the limit of any converging subsequence, and thus of $(\pi_k)_{k\in \N}$, is $\pi$. 

Finally, since $\pi_k$ has finite support and every point in its support has pairwise distinct coordinates, we may choose $0<\delta_k\leq 1/k$ such that $\pi_k(D_{\delta_k}) = 0$.
\end{proof}

\begin{rmk}\label{rmk:finite-cost-general}
\Cref{prop:finite-cost-separation} shows that finite cost alone yields finitely supported approximations which give no mass to $D_{\delta_k}$, but $\delta_k\to0$, so that no positive exclusion radius is obtained for $\pi$ itself. For minimizers, under a condition on $\kappa_\rho$, \Cref{thm:diagonal-estimate} gives a positive exclusion radius, which by \Cref{cor:diagonal-uniform} is independent of $n$. We also notice that the proof of \Cref{prop:finite-cost-separation} only relies on \Cref{thm:tchakaloff}, and therefore applies to any finite family of moment functions in $C_b(X^M)$; this will be used in \Cref{sec:semiclassical-gamma}.
\end{rmk}

Assume henceforth that $I(\rho)<+\infty$.

\begin{thr}\label{thm:diagonal-estimate} If $\kappa_\rho(\alpha) < 1/(2(M-1))$ and $\pi^n \in \Pi_A^n(\rho)$ is optimal in \eqref{eq:repulsive-moment-problem} then, \[\pi^n(D_{\alpha_n}) = 0\]where $\alpha_n>0$ is such that \[h(\alpha_n)= {\dfrac{I(\rho) - I^n(\rho)}{1-2(M-1)\kappa_\rho(\alpha)} + 2(M-1)h(\alpha)}.\]
\end{thr}

\begin{proof}
The definition of $\alpha_n$ and the strict decrease of $h$ give $0<\alpha_n<\alpha$. Suppose initially that $\pi^n$ is finitely supported. Assume that $\pi^n(D_{\alpha_n}) >0$. We will show that this allows us to construct a better competitor for \eqref{eq:repulsive-moment-problem}, contradicting the optimality of $\pi^n$.
Let $\bar x=(\bar x_1,\dots,\bar x_M)\in\supp(\pi^n)\cap D_{\alpha_n}$ and let $a=\pi^n(\{\bar x\})$. After relabeling the coordinates, we may assume that $d(\bar x_1,\bar x_2)<\alpha_n$. Let $\pi^\star$ be optimal for \eqref{eq:repulsive-problem} and define the ``good'' set
\[
    \mathcal G_{\bar x}
    =\Bigl\{(z_1,\dots,z_M)\in X^M:
    d(z_1,\bar x_j)\geq\alpha,\quad
    d(\bar x_1,z_j)\geq\alpha
    \text{ for }j=2,\dots,M\Bigr\}.
\]
We have the following estimate for the mass of $\mathcal G_{\bar x}$:
\[
\begin{aligned}
    \pi^\star\bigl(X^M\setminus\mathcal G_{\bar x}\bigr)
    &=\pi^\star\Biggl(
    \bigcup_{j=2}^M\Bigl[
    \bigl(B(\bar x_j,\alpha)\times X^{M-1}\bigr)
    \cup\bigl(X^{j-1}\times B(\bar x_1,\alpha)\times X^{M-j}\bigr)
    \Bigr]\Biggr)\\
    &\leq\sum_{j=2}^M\Bigl[
    \rho\bigl(B(\bar x_j,\alpha)\bigr)
    +\rho\bigl(B(\bar x_1,\alpha)\bigr)\Bigr]\\
    &\leq2(M-1)\kappa_\rho(\alpha).
\end{aligned}
\]
Hence, from $\kappa_\rho(\alpha)<1/(2(M-1))$,
\begin{equation}\label{eq:diagonal-good-mass}
    \pi^\star(\mathcal G_{\bar x})
    \geq1-2(M-1)\kappa_\rho(\alpha)>0.
\end{equation}
Define the following transport plan:
\begin{equation}\label{eq:diagonal-mixture}
    \pi_\epsilon = (1-\epsilon) {\pi}^n + \epsilon \pi^\star, \qquad \epsilon \in (0,1).
\end{equation}
Clearly, $\pi_\epsilon \in \Pi_A^n(\rho)$. Moreover, we can explicitly estimate its cost:
\begin{equation}\label{eq:diagonal-mixture-cost}
    \begin{aligned}
        \int_{X^M} c \d \pi_\epsilon &= (1-\epsilon) I^n(\rho) + \epsilon I(\rho) 
        =I^n(\rho) + \epsilon (I(\rho) - I^n(\rho)).
    \end{aligned}
\end{equation}

Since $d(\bar x_1,\bar x_2)<\alpha_n<\alpha$, we have $(\bar x_1,\dots, \bar x_M)\notin\mathcal G_{\bar x}$. Choose $\epsilon \in (0,1)$ such that $\epsilon \pi^\star(\mathcal G_{\bar x}) \leq (1-\epsilon) a$  and define \[\tilde \pi_\epsilon = \pi_\epsilon  - \epsilon \pi^\star(\mathcal G_{\bar x}) \delta_{(\bar x_1,\dots,\bar x_M)} - \epsilon \pi^\star\vert_{\mathcal G_{\bar x}} + \epsilon (T_1)_\#\bigl(\pi^\star\vert_{\mathcal G_{\bar x}}\bigr) + \epsilon (T_2)_{\#}\bigl(\pi^\star\vert_{\mathcal G_{\bar x}}\bigr),\] where $T_1,T_2: \mathcal G_{\bar x} \to X^M$ are given by $T_1(z_1,\dots,z_M) = (\bar x_1,z_2,\dots,z_M)$ and $T_2(z_1,\dots,z_M) = (z_1,\bar x_2,\dots, \bar x_M)$.

Observe that $\tilde \pi_\epsilon$ is a nonnegative measure, since $\epsilon \pi^\star(\mathcal G_{\bar x})\delta_{(\bar x_1,\dots,\bar x_M)} + \epsilon \bigl(\pi^\star\vert_{\mathcal G_{\bar x}}\bigr)$ is a submeasure of $\pi_\epsilon$ by our choice of $\epsilon$. Moreover, $\tilde \pi_\epsilon$ is a probability measure. We also have

\begin{equation}
    \begin{aligned}
        (p_1)_{\#}\tilde \pi _\epsilon &= (p_1)_{\#}\pi_\epsilon - \epsilon \pi^\star(\mathcal G_{\bar x}) \delta_{\bar x_1} - \epsilon (p_1)_{\#}\bigl(\pi^\star\vert_{\mathcal G_{\bar x}}\bigr) \\  &\quad +\epsilon (p_1)_{\#} (T_1)_{\#}\bigl(\pi^\star\vert_{\mathcal G_{\bar x}}\bigr) + \epsilon (p_1)_{\#}(T_2)_{\#}\bigl(\pi^\star\vert_{\mathcal G_{\bar x}}\bigr) \\
        &= (p_1)_{\#}\pi_\epsilon - \epsilon \pi^\star(\mathcal G_{\bar x}) \delta_{\bar x_1} - \epsilon(p_1)_{\#}\bigl(\pi^\star\vert_{\mathcal G_{\bar x}}\bigr) \\  &\quad +\epsilon \pi^\star(\mathcal G_{\bar x})\delta_{\bar x_1} + \epsilon (p_1)_{\#} \bigl(\pi^\star\vert_{\mathcal G_{\bar x}}\bigr) \\
        &=(p_1)_{\#} \pi_\epsilon.
    \end{aligned}
\end{equation}
and similarly $(p_j)_{\#} \tilde \pi_{\epsilon} = (p_j)_{\#}\pi_\epsilon$ for $j=2,\dots,M$. Since $\tilde \pi_\epsilon$ and $\pi_\epsilon$ have the same marginals, they have the same moments against $\varphi_i\circ p_j$ for every $i=1,\dots,n$ and $j=1,\dots, M$. Moreover, since $\Theta(x_1,\dots, x_M)=\theta(x_1)+\dots +\theta(x_M)$, we have
\[\int_{X^M}\Theta \d\tilde \pi_\epsilon=\int_{X^M}\Theta \d\pi_\epsilon\leq A.\]
Therefore $\tilde \pi_\epsilon \in \Pi_A^n(\rho)$.
Integrating the cost $c$ against the measure $\tilde \pi_\epsilon$ we obtain

\begin{equation}\label{eq:diagonal-swap-cost}
    \begin{aligned}
        \int_{X^M} &c \d \tilde \pi_\epsilon \\ &= \int_{X^M} c\d\pi_{\epsilon} - \epsilon \pi^\star(\mathcal G_{\bar x})\int_{X^M} c \d \delta_{(\bar x_1,\dots , \bar x_M)}   -\epsilon \int_{\mathcal G_{\bar x}} c \d\pi^\star \\& \quad + \epsilon \int_{\mathcal G_{\bar x}} c\circ T_1\d\pi^\star + \epsilon \int_{\mathcal G_{\bar x}} c\circ T_2 \d\pi^\star \\
        &=\int_{X^M} c \d\pi_{\epsilon}  - \epsilon \int_{\mathcal G_{\bar x}} c(\bar x_1,\dots,\bar x_M) \d\pi^\star  - \epsilon \int_{\mathcal G_{\bar x}} c(z_1,\dots,z_M) \d\pi^\star \\ & \quad + \epsilon \int_{\mathcal G_{\bar x}} c(\bar x_1,z_2,\dots,z_M) \d\pi^\star + \epsilon \int_{\mathcal G_{\bar x}}c(z_1,\bar x_2,\dots,\bar x_M) \d \pi^\star  \\
        &= \int_{X^M} c \d\pi_\epsilon +\epsilon \sum^M_{j=2} \int_{\mathcal G_{\bar x}} [ -h(d(\bar x_1,\bar x_j)) - h(d(z_1,z_j)) \\ & \hspace{11em} + h(d(\bar x_1,z_j)) +h(d(z_1,\bar x_j)) ]\d\pi^\star.
    \end{aligned}
\end{equation}

Since $d(\bar x_1,\bar x_2)<\alpha_n$ and for $(z_1,\dots,z_M) \in \mathcal G_{\bar x}$ we have $d(z_1,\bar x_j) \geq \alpha$ and $d(\bar x_1,z_j) \geq \alpha$ for every $j=2,\dots,M$, we can estimate the sum of the integrands as
\begin{equation}\label{eq:diagonal-cost-bound}\begin{aligned}
    \sum^M_{j=2} \Bigl[ -h(d(\bar x_1,\bar x_j)) - h(d(z_1,z_j)) &+ h(d(\bar x_1,z_j)) +h(d(z_1,\bar x_j)) \Bigr] \\
    &\leq -h(d(\bar x_1,\bar x_2)) + 2(M-1)h(\alpha) \\
    &< -h(\alpha_n) + 2(M-1)h(\alpha).
    \end{aligned}
\end{equation}

Substituting \eqref{eq:diagonal-cost-bound} and \eqref{eq:diagonal-mixture-cost} into \eqref{eq:diagonal-swap-cost} we obtain

\begin{equation}\begin{aligned}
    \int_{X^M} c \d \tilde \pi_\epsilon &< \int_{X^M} c \d \pi_\epsilon + \epsilon \pi^\star(\mathcal G_{\bar x})\left( -h(\alpha_n) +2(M-1)h(\alpha) \right) \\
    &\leq I^n(\rho) + \epsilon\left[ I(\rho) - I^n(\rho) +  \pi^\star(\mathcal G_{\bar x}) \left( -h(\alpha_n) +2(M-1)h(\alpha)\right) \right].
    \end{aligned}
\end{equation} Since, by assumption, $h(\alpha_n) = \dfrac{I(\rho) - I^n(\rho)}{1-2(M-1)\kappa_\rho(\alpha)} + 2(M-1)h(\alpha)$, combining this with \eqref{eq:diagonal-good-mass} yields

\begin{equation}
    \begin{aligned}
        \int_{X^M} c \d\tilde \pi_\epsilon &< I^n(\rho) + \epsilon \left[ I(\rho) - I^n(\rho) - \pi^\star(\mathcal G_{\bar x}) \left( \frac{I(\rho) - I^n(\rho)}{1-2(M-1)\kappa_\rho(\alpha)} \right)\right] \\
        & \leq I^n(\rho) + \epsilon\left[ I(\rho) - I^n(\rho) - (1-2(M-1)\kappa_\rho(\alpha)) \left(\frac{I(\rho) - I^n(\rho)}{1-2(M-1)\kappa_\rho(\alpha)}\right)\right] \\
        & = I^n(\rho)
    \end{aligned}
\end{equation}This is a contradiction, since $\tilde \pi_\epsilon \in \Pi_A^n(\rho)$.

Let now $\pi^n$ be an arbitrary minimizer. By \Cref{thm:extreme-minimizer}, there exists a sequence of finitely supported minimizers $(\pi^{n;k})_{k\in \N}$ converging narrowly to $\pi^n$. The preceding argument gives us that $\pi^{n;k}(D_{\alpha_n}) = 0 $ for every $k \in \N$. Since $D_{\alpha_n}$ is open and does not depend on $k$, we have \[0\leq \pi^n(D_{\alpha_n})\leq \liminf_{k\to \infty} \pi^{n;k}(D_{\alpha_n}) = 0.\]This concludes the proof.
\end{proof}

\begin{rmk}If $\bigcup_n\I^n$ fully characterizes $\rho$ and the admissible sets $\Pi^n_A$ satisfy the limit condition in \Cref{thm:moment-gamma}, then $I^n(\rho)\to I(\rho)$. Consequently, the exclusion radii in \Cref{thm:diagonal-estimate} satisfy $h(\alpha_n)\to2(M-1)h(\alpha)$. Thus, in the limit one recovers the diagonal exclusion radius (see \cite[Theorem 1.3 (i)]{ColoDiMaStra2019SIMA}).\end{rmk}




Since $0\leq I(\rho)-I^n(\rho)\leq I(\rho)$ and $h$ is decreasing, the diagonal exclusion is uniform in $n$:
\begin{cor}\label{cor:diagonal-uniform}
If $\kappa_\rho(\alpha)<1/(2(M-1))$ for some $\alpha>0$, set
\[
    \alpha^\star=h^{-1}\left(
    \frac{I(\rho)}{1-2(M-1)\kappa_\rho(\alpha)}+2(M-1)h(\alpha)\right).
\]
Then, for every $n\in\N$, every minimizer $\pi^n$ of \eqref{eq:repulsive-moment-problem} satisfies $\pi^n(D_{\alpha^\star})=0$.
\end{cor}
\noindent For $0<\bar\alpha<\alpha^\star$, define $\bar h(r) = \min(h(r),h(\bar \alpha))$, the truncated cost function $\bar c$ defined as in \eqref{eq:repulsive-cost} with respect to the function $\bar h$, and $\bar I^n(\rho)$ and $\bar I(\rho)$ defined as in \eqref{eq:repulsive-moment-problem} and \eqref{eq:repulsive-problem}, respectively, with respect to $\bar c$. 

\begin{thr}\label{thm:repulsive-truncation} $\pi^n\in \Pi^n_A(\rho)$ is optimal for $I^n(\rho)$ if, and only if, it is optimal for $\bar I^n(\rho)$.
\end{thr}

\begin{proof} We first show that every minimizer $\bar\pi^n$ of $\bar I^n(\rho)$ satisfies $\bar\pi^n(D_{\alpha^\star})=0$. As in the proof of \Cref{thm:diagonal-estimate}, one may initially assume that $\bar \pi^n$ is finitely supported and then pass to the limit. The only change concerns \eqref{eq:diagonal-cost-bound}: Since $d(\bar x_1,\bar x_2) < \alpha^\star$ and $\bar \alpha < \alpha^\star$ we have  $\bar h(d(\bar x_1,\bar x_2)) > h(\alpha^\star)$, while on $\mathcal{G}_{\bar x}$, both $\bar h(d(\bar x_1,z_j))$ and $\bar h(d(z_1,\bar x_j))$ are at most $h(\alpha)$ for every $j=2,\dots,M$. Since $0\leq\bar I(\rho)-\bar I^n(\rho)\leq I(\rho)$, the definition of $\alpha^\star$ then yields the same contradiction.


Let $\pi^n$ and $\bar\pi^n$ be optimal for $I^n(\rho)$ and $\bar I^n(\rho)$, respectively. Since $\bar c\leq c$ pointwise and $\bar\pi^n(D_{\alpha^\star})=0$, we obtain
\[
\begin{aligned}
    \int_{X^M}\bar c\,\d\bar\pi^n
    &\leq\int_{X^M}\bar c\,\d\pi^n
    \leq\int_{X^M}c\,\d\pi^n
    \leq\int_{X^M}c\,\d\bar\pi^n
    =\int_{X^M}\bar c\,\d\bar\pi^n,
\end{aligned}
\]
hence all the inequalities are equalities, and the result follows.
\end{proof}

The equivalence obtained in \Cref{thm:repulsive-truncation} has been used in the literature to obtain bounds on the Kantorovich potential of the corresponding dual problem (see \cite{ColoDiMaStra2019SIMA, ButtChamDePa2018AMO}). For the moment-constrained problems considered in this section, a similar result is possible but requires stronger assumptions:
\begin{assumption}\label{ass:compact-local}
The space $X$ is compact and $\rho \in\P(X)$ has full support on $X$. For fixed $n$, the moment functions $\I^n=\{\varphi_1,\dots,\varphi_n\}\subset C(X)$ satisfy \eqref{eq:local-moments}, with the compact set taken to be $X$ and nodes $z_1,\dots,z_n$, and $\varphi_i(z_j)=\delta_{ij}$ for every $i,j=1,\dots,n$.
\end{assumption}
\noindent The compactness of $X$ permits us to choose $\Theta \equiv 0$. Set $E_n=\operatorname{span}(\I^n)$ and define $U_i\coloneq\{\varphi_i>0\},$ $ J_i\coloneq \{j\in\{1,\dots,n\}:U_i\cap U_j\neq\emptyset\},$ $$R_i\coloneq
    \frac{\displaystyle\sum_{j\in J_i}
        \int_X\varphi_j\,\d\rho}
    {\displaystyle\int_X\varphi_i\,\d\rho},$$ 
which measures the total $\rho$-mass assigned to the moment functions whose positivity set intersects the positivity set of $\varphi_i$ relative to the mass assigned to $\varphi_i$ itself. The denominator is positive because $\rho$ has full support and $\varphi_i(z_i)=1$. Since $\bar c$ is continuous on the compact space $X^M$, the function
\[
    \omega(t)\coloneq
    \sup_{\substack{x_1,x_1',x_2,\dots,x_M \in X \\ d(x_1,x_1')\leq t} }
        |\bar c(x_1,x_2,\dots,x_M)-\bar c(x_1',x_2,\dots,x_M)|
\]
satisfies $\omega(t)\to0$ as $t\to0$.

\begin{thr}\label{thm:repulsive-potentials} Under \Cref{ass:compact-local}, the dual problem to $I^n(\rho)$, namely
\begin{equation}\label{eq:repulsive-dual}
    \sup_{\substack{
        u_1,\dots,u_M\in E_n\\
        \sum_{j=1}^M u_j(x_j)\leq c(x_1,\dots,x_M)
    }}
    \left\{
        \sum_{j=1}^M\int_X u_j\,\d\rho
    \right\}
\end{equation}admits a maximizer of the form $(u,\dots,u)$ where
$u=\sum_{i=1}^n a_i\varphi_i$ and
\[
    -\frac{\| \bar c \|_\infty}{M}-\omega(r_n)R_i
    \leq a_i\leq\frac{\| \bar c \|_\infty}{M},
    \qquad i=1,\dots,n.
\]
\end{thr}

\begin{proof}
We only treat the case $M=2$; the case $M>2$ is analogous. Since $\bar c$ is bounded
and continuous and $\rho\otimes\rho$ has full support,
\Cref{thm:ot-duality} gives attainment in the dual problem
with cost $\bar c$. By symmetry of the cost, if $(u_1,u_2)$
is a maximizer, then $(u,u)$ with $u=(u_1+u_2)/2$ is also
a maximizer.

Let $(u,u)$ be such a maximizer, with
$u=\sum_{i=1}^n a_i\varphi_i$. Taking $x=y$ in the dual
constraint gives
$2u(x)\leq\bar c(x,x)=\|\bar c\|_\infty$.
Since $u(z_i)=a_i$, we obtain
$a_i\leq\|\bar c\|_\infty/2$ for every $i$.

To prove the lower bound, fix $i\in\{1,\dots,n\}$.
If $a_i\geq-\|\bar c\|_\infty/2$, there is nothing to prove.
Otherwise, define
\[
    \widetilde u
    =u+\left(-\frac{\|\bar c\|_\infty}{2}-a_i\right)\varphi_i.
\]
This replaces the coefficient $a_i$ by
$-\|\bar c\|_\infty/2$ and leaves the other coefficients
unchanged. In particular,
$\widetilde u(z_i)=-\|\bar c\|_\infty/2$, whereas
$\widetilde u(z_j)=u(z_j)$ for $j\neq i$.

Let us estimate the error in the admissibility of
$(\widetilde u,u)$. For $j\neq i$, feasibility of $(u,u)$
gives
$\widetilde u(z_j)+u(y)\leq\bar c(z_j,y)$.
The same inequality holds for $j=i$, since
$\widetilde u(z_i)+u(y)=-\|\bar c\|_\infty/2+u(y)\leq0$
and $\bar c$ is nonnegative.
Multiplying these inequalities by $\varphi_j(x)$ and
summing over $j$, we obtain
\[
    \widetilde u(x)+u(y)
    \leq\sum_{j=1}^n\bar c(z_j,y)\varphi_j(x)
    \leq\bar c(x,y)+\omega(r_n),
\]
where the last inequality follows from $\supp(\varphi_j) \subset B(z_j,r_n)$. Moreover, $\widetilde u=u$ outside
$U_i$, so the error is confined to $U_i$:
\begin{equation}\label{eq:repulsive-dual-error}
    \widetilde u(x)+u(y)
    \leq\bar c(x,y)+\omega(r_n)\chi_{U_i}(x).
\end{equation}

We correct this error by setting
$\zeta_i=\sum_{j\in J_i}\varphi_j$.
The partition of unity property gives $0\leq\zeta_i\leq1$
and $\zeta_i=1$ on $U_i$. Indeed, whenever $x\in U_i$ and
$\varphi_j(x)>0$, we have $j\in J_i$.
Thus $\chi_{U_i}\leq\zeta_i$, and
\eqref{eq:repulsive-dual-error} shows that
$(\widetilde u-\omega(r_n)\zeta_i,u)$ is admissible for
the dual problem with cost $\bar c$.

Comparing its value with that of the maximizer $(u,u)$ gives
\[
\begin{aligned}
    0
    &\geq
    \int_X\bigl(\widetilde u-\omega(r_n)\zeta_i-u\bigr)\,\d\rho\\
    &=
    \left(-\frac{\|\bar c\|_\infty}{2}-a_i\right)
    \int_X\varphi_i\,\d\rho
    -\omega(r_n)
    \sum_{j\in J_i}\int_X\varphi_j\,\d\rho.
\end{aligned}
\]
Dividing by $\int_X\varphi_i\,\d\rho>0$ yields
$a_i\geq-\|\bar c\|_\infty/2-\omega(r_n)R_i$,
as required.

Finally, since $\bar c\leq c$, the pair $(u,u)$ is also
admissible for the original dual problem.
By \Cref{thm:repulsive-truncation}, its value is
$2\int_Xu\,\d\rho=\bar I^n(\rho)=I^n(\rho)$.
Weak duality therefore shows that $(u,u)$ is optimal for
the original cost $c$.
\end{proof}

\begin{rmk}
Since the functions $\varphi_i$ are nonnegative and sum
to one, the coefficient bounds imply
\[
    \|u\|_\infty
    \leq\frac{\|\bar c\|_\infty}{M}
    +\omega(r_n)\max_{1\leq i\leq n}R_i.
\]
For the Coulomb cost,
$\|\bar c\|_\infty=\binom{M}{2}/\bar\alpha$ and $\bar c$ is
$(M-1)\bar\alpha^{-2}$-Lipschitz in each variable. Consequently,
\[
    \|u\|_\infty
    \leq\frac{M-1}{2\bar\alpha}
    +\frac{(M-1)r_n}{\bar\alpha^2}
    \max_{1\leq i\leq n}R_i.
\]
In particular, for fixed $M$ and $\bar\alpha$, if $r_n\to0$
and the ratios $R_i$ are bounded uniformly in $n$ and $i$,
the error term tends to zero.
\end{rmk}

\begingroup
\ExplSyntaxOn
\cs_set_eq:NN \span \tex_span:D
\ExplSyntaxOff
\section{A semiclassical \texorpdfstring{$\Gamma$}{Gamma}-convergence result}
\label{sec:semiclassical-gamma}

We now explain how the moment-constrained Coulomb problem of \Cref{sec:repulsive-costs} and, more generally, transport problems with finitely many moment equalities arise as the
zero-order limit of a fermionic Lieb functional.  We restrict the discussion to the
finite-dimensional moment problem, that is, to one fixed family $W^n$\corr{ of moment functions encoding the equality constraints, as specified below}.
We consider\corr{ measures of mass one satisfying the confinement bound}.  In the following, we present the
semiclassical analogue of the regularization results for the Lieb and Hohenberg--Kohn
functionals in \cite{BindDePa2017JEP,Lewi2018CRM,CotarFriKlu2018}, combined with the
finite-moment transport framework of \cite{AlfCoyEhrLom2021Approx}.

Let $M\geq2$, $Z=(\mathbb R^3)^M$, and let $\mathcal P_{\mathrm{sym}}(Z)$ be the
narrowly closed set of probability measures invariant under permutations of the particle
coordinates.  We write $\operatorname{Sym}\pi=(M!)^{-1}\sum_{\sigma\in S_M}\sigma_\#\pi$.
The cost $c$ of the preceding sections is now the Coulomb cost
\[
c(\corr{\mathbf{x}})=\sum_{1\leq i<j\leq M}\frac{1}{|x_i-x_j|},
 \qquad \corr{\mathbf{x}}=(x_1,\ldots,x_M),
\]
with the convention $c=+\infty$ whenever two particle coordinates coincide.  We assume
that $\Theta:Z\to[0,+\infty)$ is continuous, coercive and permutation invariant, and fix
$A>0$.

For the equality constraints, let $g_1,\ldots,g_q\in C_b(Z)$
be permutation-invariant moment functions and let
$b\in\R^q$ be the prescribed vector.  We write
\[
\Psi_n(\pi):=\left(\int_Zg_1\d\pi,\ldots,\int_Zg_q\d\pi\right),
 {\qquad
 W^n:=\{g_i-b_i,\ b_i-g_i:\ i=1,\ldots,q\}.}
\]
Here $n$ labels the moment family and $q=q_n$ counts the moment equalities.
With the moment functions $W^n$ and mass $m=1$, the inequalities in \eqref{eq:moment-set} encode precisely $\Psi_n(\pi)=b$.
In this section the admissible set includes symmetry:
\[
\Pi_{A,\mathrm{sym}}^n:=\Pi_A^n \cap\P_{\mathrm{sym}}(Z)
=\left\{\pi\in\P_{\mathrm{sym}}(Z):
       \Psi_n(\pi)=b,\ \int_Z\Theta\d\pi\leq A\right\}.
\]
We define, on $\P(Z)$,
\[
\C_{\mathrm{sym}}^n(\pi):=
 \begin{cases}
 \displaystyle\int_Z c\d\pi,&\pi\in\Pi_{A,\mathrm{sym}}^n,\\[1mm]
 +\infty,&\pi\notin{\Pi_{A,\mathrm{sym}}^n}.
 \end{cases}
\]
Thus $\C_{\mathrm{sym}}^n$ agrees with the extended functional $\C^n$ of
\Cref{sec:gamma-convergence} on $\P_{\mathrm{sym}}(Z)$, and is infinite outside this set.
When the moment equalities are prescribed by a common marginal $\rho$, that is, $g_i=M^{-1}\sum_{j=1}^M\varphi_i\circ p_j$ and $b_i=\int_{\R^3}\varphi_i\,\d\rho$, and $\Theta$ is as in \Cref{sec:repulsive-costs}, symmetrization shows that
\[
 \inf_{\pi\in\Pi_{A,\mathrm{sym}}^n}\C(\pi)=I^n(\rho).
\]
This is the symmetric $M$-particle formulation of the $M$-marginal value considered in \Cref{sec:repulsive-costs}.

Let $\mathfrak S_a$ denote the set of mixed fermionic states on
$L^2_a((\R^3)^M)$, let $\pi_\Gamma$ be the associated symmetric $M$-particle law,
and let
\[
 K(\Gamma):=\frac12\operatorname{Tr}(-\Delta\Gamma)
 \in[0,+\infty]
\]
be the kinetic energy.
For a probability measure $\pi$ on $Z$, define the partially
minimized kinetic energy
\[
 T(\pi):=\inf\big\{K(\Gamma):\Gamma\in\mathfrak S_a,
                 \ \pi_\Gamma=\pi\big\},
 \qquad \inf\varnothing:=+\infty.
\]
The lifted semiclassical functionals are
\begin{equation}
 \label{eq:lifted-semiclassical-functional}
 F_n^\varepsilon(\pi):=\varepsilon T(\pi)+\C_{\mathrm{sym}}^n(\pi),
 \qquad \varepsilon>0,
\end{equation}
and their zero-order candidate is $F_n^0:=\C_{\mathrm{sym}}^n$.  Partial minimization gives the identity
\[
 \inf_{\pi\in\P(Z)}F_n^\varepsilon(\pi)
 =\inf_{\substack{\Gamma\in\mathfrak S_a\\ \pi_\Gamma\in{\Pi_{A,\mathrm{sym}}^n}}}
 \left\{\varepsilon K(\Gamma)+\int_Z c\d\pi_\Gamma\right\}
 =:\mathcal L_n^\varepsilon.
\]

Set
\[
 H_n:=\operatorname{aff}\{\Psi_n(\pi):\pi\in\P_{\mathrm{sym}}(Z)\}.
\]
The recovery construction below replaces the atoms of a transport plan by fermionic wave packets, which slightly perturbs the moments. To recover the constraint $\Psi_n(\pi)=b$ exactly, we mix the resulting state with finitely many fixed states whose moment vectors have a convex hull containing a neighborhood of $b$ in $H_n$ (these moment vectors need not be close to $b$); since only a small mixing weight is needed, the energies are almost unchanged. The following assumption guarantees that such states exist. It is the analogue of Slater's condition in convex optimization: $b$ must lie in the relative interior, within the affine hull $H_n$, of the moment vectors of admissible finite-energy states since, on the boundary, small errors in the moments could not be corrected in every direction.

\begin{assumption}
\label{ass:moment-correctability}
There are finitely many states $\Gamma^0,\ldots,\Gamma^J\in\mathfrak S_a$ such that
\[
 K(\Gamma^j)<+\infty,\qquad
 \int_Z c\d\pi_{\Gamma^j}<+\infty,\qquad
 \int_Z\Theta \d\pi_{\Gamma^j}<A
 \quad (j=0,\ldots,J).
\]
Moreover, $b\in H_n$ and, for some $r>0$,
\[
 B(b,r)\cap H_n\subset
 \operatorname{conv}\left\{\Psi_n(\pi_{\Gamma^0}),\ldots,
                         \Psi_n(\pi_{\Gamma^J})\right\}.
\]
\end{assumption}

\begin{rmk}\label{rmk:full-support}

\Cref{ass:moment-correctability} requires $b$ to lie in the relative interior, within $H_n$, of the moment vectors of strictly confined fermionic states with finite energies, so that small errors in the moments can be corrected by convex combinations. It plays the same role as the full-support assumption of \Cref{thm:duality}. Similarly, it holds as soon as there is a state $\Gamma^\circ\in\mathfrak S_a$ with finite kinetic and Coulomb energies such that
\[
 \int_Z\Theta\d\pi_{\Gamma^\circ}<A,\qquad
 \Psi_n(\pi_{\Gamma^\circ})=b,\qquad
 \supp(\pi_{\Gamma^\circ})=Z:
\]
localizing $\Gamma^\circ$ with smooth, permutation-invariant cutoff functions produces finite-energy fermionic states whose moment vectors surround $b$.

\end{rmk}

In particular, there are weights $\lambda_j^*\geq0$, summing to one, such that
\[
 \Gamma^{\mathrm{ref}}:=\sum_{j=0}^J\lambda_j^*\Gamma^j,
 \qquad \sigma:=\pi_{\Gamma^{\mathrm{ref}}},\qquad
 \Psi_n(\sigma)=b,\quad \int_Z\Theta\d\sigma<A.
\]
Both energies of $\Gamma^{\mathrm{ref}}$ are finite.  We use this reference state below.

We next give the recovery construction for the $\Gamma$-convergence result.  We use
\[
 D:=\{\corr{\mathbf{x}}\in Z:\min_{i\neq j}|x_i-x_j|=0\},
 \qquad
 D_\delta:=\{\corr{\mathbf{x}}\in Z:\min_{i\neq j}|x_i-x_j|<\delta\},
\]
and the bounded-Lipschitz distance $d_{\mathrm{BL}}$, 
\corr{
\[
 d_{\mathrm{BL}}(\pi,\pi'):=\sup\left\{\int_Zf\d\pi-\int_Zf\d\pi':\ \|f\|_\infty\leq1,\ \operatorname{Lip}(f)\leq1\right\},
 \qquad \pi,\pi'\in\P(Z),
\]}

which metrizes narrow convergence on
$\P(Z)$.

\begin{lem}[Sparse approximation in energy]
\label{lem:sparse-section}
Let $\pi\in\Pi_{A,\mathrm{sym}}^n$ satisfy $\int_Zc\d\pi<+\infty$ and
$\int_Z\Theta\d\pi<A$.  Then there are finite symmetric atomic laws
\[
 \pi^k=\sum_{\ell=1}^{L_k}\alpha_{k,\ell}\,\operatorname{Sym}\delta_{\corr{\mathbf{x}_{k,\ell}}}
\]
with all $\corr{\mathbf{x}_{k,\ell}}\notin D$, such that
\[
 d_{\mathrm{BL}}(\pi^k,\pi)\to0,\qquad
 \int_Z c\d\pi^k=\int_Z c \d\pi,\qquad
 \int_Z\Theta \d\pi^k=\int_Z\Theta \d\pi,
 \qquad \Psi_n(\pi^k)=b.
\]
\end{lem}
\begin{proof}

By \Cref{rmk:finite-cost-general}, the proof of \Cref{prop:finite-cost-separation} applies with the moment functions $g_1,\dots,g_q$ in place of $\varphi_i\circ p_j$. It gives finitely supported measures $\widehat\pi^k$, converging narrowly to $\pi$, with atoms in $\supp(\pi)\setminus D$ and with the same integrals as $\pi$ against $g_1,\dots,g_q$, $c$ and $\Theta$. Since $\pi$ is symmetric and $g_1,\dots,g_q$, $c$ and $\Theta$ are permutation invariant, the symmetrized measures $\pi^k\coloneq\operatorname{Sym}\widehat\pi^k$ have the same properties, converge narrowly to $\operatorname{Sym}\pi=\pi$, and are of the stated form. Finally, $d_{\mathrm{BL}}$ metrizes narrow convergence.

\end{proof}

\begin{lem}[Packet representation of a sparse law]
\label{lem:packets-section}
Let $\pi^k=\sum_{\ell=1}^{L_k}\alpha_{k,\ell}\operatorname{Sym}\delta_{\corr{\mathbf{x}_{k,\ell}}}$ be as in
Lemma~\ref{lem:sparse-section} and set
\[
 d_k:=\min_{\ell}\min_{i\neq j}|x_{k,\ell,i}-x_{k,\ell,j}|>0.
\]
Fix a real, even $\chi\in C_c^\infty(B(0,1))$ with $\|\chi\|_2=1$ and write
$\chi_{\eta,a}(x)=\eta^{-3/2}\chi((x-a)/\eta)$.  For $0<\eta<d_k/3$ there is a mixed
fermionic state $\Gamma_{k,\eta}$ such that
\begin{align*}
 K(\Gamma_{k,\eta})&=\frac{M}{2\eta^2}\int_{\R^3}|\nabla \chi|^2,\\
 \pi_{\Gamma_{k,\eta}}&=\sum_{\ell=1}^{L_k}\alpha_{k,\ell}\operatorname{Sym}
             \bigotimes_{i=1}^M\chi_{\eta,x_{k,\ell,i}}^2.
\end{align*}
Consequently, as $\eta\downarrow0$,
\[
 \pi_{\Gamma_{k,\eta}}\weakto\pi^k,
 \qquad
 \int_Zc\d\pi_{\Gamma_{k,\eta}}\to\int_Z c\d\pi^k,
 \qquad
 \Psi_n(\pi_{\Gamma_{k,\eta}})\to\Psi_n(\pi^k).
\]
If $\int_Z\Theta \d\pi^k<A$, then the confinement inequality also holds for $\eta$ small.
\end{lem}
\begin{proof}
For each $\ell$, the supports of the orbitals $\chi_{\eta,x_{k,\ell,i}}$ are disjoint, hence
they are orthonormal.  Let $\mathcal S_{k,\ell,\eta}$ be the associated Slater determinant and
set
\[
 \Gamma_{k,\eta}:=\sum_{\ell=1}^{L_k}\alpha_{k,\ell}
 |\mathcal S_{k,\ell,\eta}\rangle\langle\mathcal S_{k,\ell,\eta}|.
\]
Disjointness removes all exchange terms and gives the displayed formula for the law.  The
kinetic identity follows by scaling (see \cite{Lewi2018CRM} for a similar construction).  The remaining convergences follow from continuity of
$c$ away from $D$, bounded continuity of the moments, and continuity of $\Theta$ on the
compact support.
\end{proof}

\begin{lem}[Exact finite-dimensional correction]
\label{lem:correction-section}
Let $(\Gamma_\eta)$ be states with $K(\Gamma_\eta)\le A_0\eta^{-2}$,
$\int_Zc\d\pi_{\Gamma_\eta}\le B_0$, $\int_Z\Theta \d\pi_{\Gamma_\eta}<A$, and
\[
 \beta_\eta:=|\Psi_n(\pi_{\Gamma_\eta})-b|\longrightarrow0.
\]
Under Assumption~\ref{ass:moment-correctability}, there are exactly feasible states
$\widetilde\Gamma_\eta$ satisfying
\[
 \Psi_n(\pi_{\widetilde\Gamma_\eta})=b,
 \qquad
 K(\widetilde\Gamma_\eta)\le A_0\eta^{-2}+C\beta_\eta,
 \qquad
 \left|\int_Zc\d\pi_{\widetilde\Gamma_\eta}-\int_Zc\d\pi_{\Gamma_\eta}\right|
 \le C\beta_\eta,
\]
 and $d_{\mathrm{BL}}(\pi_{\widetilde\Gamma_\eta},\pi_{\Gamma_\eta})\to0$.
\end{lem}
\begin{proof}
Let $v_j:=\Psi_n(\pi_{\Gamma^j})$.  For $s\in(0,1]$ and a probability vector
$\lambda=(\lambda_0,\ldots,\lambda_{{J}})$, set
\[
 \widetilde\Gamma_\eta(s,\lambda):=(1-s)\Gamma_\eta
       +s\sum_{j=0}^{{J}}\lambda_j\Gamma^j.
\]
If $\delta_\eta:=b-\Psi_n(\pi_{\Gamma_\eta})$, exact feasibility is equivalent to
\[
 \sum_{j=0}^{J}\lambda_jv_j
 =b+\frac{1-s}{s}\delta_\eta.
\]
If $\beta_\eta=0$, leave the state unchanged.  Otherwise,
choose $s_\eta:=\min\{1,2\beta_\eta/r\}$.  For small $\eta$ the right-hand side belongs
to $B(b,r)\cap H_n$, so the assumed convex-hull condition gives a suitable $\lambda_\eta$.
Convexity preserves the confinement inequality.  Since kinetic and Coulomb energies are
affine in the state, and $s_\eta=O(\beta_\eta)$, the two stated estimates follow.  Finally,
\[
 d_{\mathrm{BL}}(\pi_{\widetilde\Gamma_\eta},\pi_{\Gamma_\eta})\le2s_\eta\to0.
\]
\end{proof}

\begin{prop}[Density of exactly feasible finite-energy laws]
\label{prop:dense-section}
Under Assumption~\ref{ass:moment-correctability}, for every
$\pi\in\Pi_{A,\mathrm{sym}}^n$ with $\int_Zc\d\pi<+\infty$ there are exactly feasible states
$\widetilde\Gamma_j$ such that, with $\widetilde\pi_j:=\pi_{\widetilde\Gamma_j}$,
\[
 d_{\mathrm{BL}}(\widetilde\pi_j,\pi)\to0,
 \qquad
 \int_Zc\d\widetilde\pi_j\to\int_Zc\d\pi,
 \qquad K(\widetilde\Gamma_j)<+\infty.
\]
\end{prop}
\begin{proof}
Let $\sigma:=\pi_{\Gamma^{\mathrm{ref}}}$ and, for $t\in(0,1]$, set
\[
 \pi_t:=(1-t)\pi+t\sigma.
\]
Then $\Psi_n(\pi_t)=b$, $\pi_t\weakto\pi$ and
$\int c\d\pi_t\to\int c \d\pi$, while
$A-\int\Theta \d\pi_t>0$.  Put $t_j:=1/j$ and
$\sigma_j:=A-\int_Z\Theta \d\pi_{t_j}$.  For each fixed $j$, apply
Lemma~\ref{lem:sparse-section} to $\pi_{t_j}$.  Then choose one of its atomic laws, say
$\pi^{k(j)}_j$, sufficiently close in distance, energy, moments and confinement,
and write $d_{j,k(j)}$ for its separation.  Choose a packet radius
\[
 0<\eta_j<\min\{d_{j,k(j)}/3,1/j\},
\]
so small that the packet state $\Gamma_{j,\eta_j}$ satisfies
\begin{align*}
 d_{\mathrm{BL}}(\pi_{\Gamma_{j,\eta_j}},\pi_{t_j})&\le 1/j,\\
 \left|\int_Z c \d\pi_{\Gamma_{j,\eta_j}}-\int_Z c \d\pi_{t_j}\right|&\le 1/j,\\
 |\Psi_n(\pi_{\Gamma_{j,\eta_j}})-b|&\le 1/j,\\
 \int_Z\Theta\d\pi_{\Gamma_{j,\eta_j}}&\le A-\sigma_j/2.
\end{align*}
This is possible because the sparse approximations converge in energy and the packet
construction converges for each fixed atomic law.  Apply Lemma~\ref{lem:correction-section}
to $\Gamma_{j,\eta_j}$ and denote the corrected state by $\widetilde\Gamma_j$.  Convexity
preserves the confinement inequality, and the correction changes the bounded-Lipschitz
distance and Coulomb energy by $o(1)$.  Therefore
\[
 d_{\mathrm{BL}}(\pi_{\widetilde\Gamma_j},\pi_{t_j})\to0,
 \qquad
 \int_Z c \d\pi_{\widetilde\Gamma_j}-\int_Zc\d\pi_{t_j}\to0.
\]
Finally, $d_{\mathrm{BL}}(\pi_{t_j},\pi)\le2/j$ and
\[
 \left|\int_Z c\d\pi_{t_j}-\int_Zc\d\pi\right|
 \le\frac1j\left|\int_Zc\d\sigma-\int_Zc\d\pi\right|.
\]
The triangle inequality gives the two asserted convergences, while
$K(\widetilde\Gamma_j)<+\infty$ follows from the packet and correction estimates.
\end{proof}

\begin{thr}[Semiclassical $\Gamma$-convergence]
\label{thm:semiclassical-gamma}
Under the standing hypotheses of this section, including continuity of $\Theta$ and
Assumption~\ref{ass:moment-correctability}, fix $n$.  Then, as
$\varepsilon\downarrow0$, the functionals $F_n^\varepsilon$ $\Gamma$-converge to $F_n^0=\C_{\mathrm{sym}}^n$
on $\P(Z)$ endowed with the narrow topology, where $\C_{\mathrm{sym}}^n=+\infty$ outside
the symmetric admissible set $\Pi_{A,\mathrm{sym}}^n$.  Moreover, the family
$(F_n^\varepsilon)_{0<\varepsilon\leq1}$ is equicoercive.

Consequently,
\[
 \mathcal L_n^\varepsilon\longrightarrow
 \min_{\pi\in\P(Z)}{\C_{\mathrm{sym}}^n}(\pi),
\]
which equals $I^n(\rho)$ in the marginal-moment setting.
If $\Gamma_\varepsilon$ are minimizers (or $o(1)$-minimizers) of the semiclassical
problem, then the family $(\pi_{\Gamma_\varepsilon})_\varepsilon$ is narrowly precompact and
every cluster point is a minimizer of $\C_{\mathrm{sym}}^n$.
\end{thr}

\begin{proof}
We first prove the $\Gamma$-liminf inequality.  Let
$\pi_\varepsilon\weakto\pi$ and suppose that
\[
 \liminf_{\varepsilon\to0}F_n^\varepsilon(\pi_\varepsilon)<+\infty.
\]
Up to a subsequence, $\pi_\varepsilon\in\Pi_{A,\mathrm{sym}}^n$.
Symmetry passes to the narrow limit.
Since all functions in $W^n$ are bounded and continuous,
the moment equalities pass to the limit.
The mass constraint is preserved by narrow convergence, and lower semicontinuity of $\Theta$ gives
\[
 \int_Z\Theta\d\pi\leq\liminf_{\varepsilon\to0}
 \int_Z\Theta\d\pi_\varepsilon\leq A.
\]
Consequently $\pi\in\Pi_{A,\mathrm{sym}}^n$.  Since $T\geq0$ and $c$ is lower semicontinuous and
bounded from below, the Portmanteau theorem yields
\[
 \liminf_{\varepsilon\to0}F_n^\varepsilon(\pi_\varepsilon)
 \geq\liminf_{\varepsilon\to0}\int_Zc\d\pi_\varepsilon
 \geq\int_Zc\d\pi
 =F_n^0(\pi).
\]

For the $\Gamma$-limsup, let
\[
 \mathcal D_n:=\left\{\pi_\Gamma:\Gamma\in\mathfrak S_a,
 K(\Gamma)<+\infty,\ \pi_\Gamma\in\Pi_{A,\mathrm{sym}}^n,\ \int_Zc\d\pi_\Gamma<+\infty\right\}.
\]
If $Q=\pi_\Gamma\in\mathcal D_n$, the constant sequence gives
\[
 \lim_{\varepsilon\downarrow0}F_n^\varepsilon(Q)
 =\int_Zc\d Q=F_n^0(Q).
\]
By Proposition~\ref{prop:dense-section}, $\mathcal D_n$ is dense in the domain of $F_n^0$
with convergence of the Coulomb energy.  Let $F'$ denote the $\Gamma$-upper limit.  Since
$F'$ is lower semicontinuous, for every $\pi$ with $F_n^0(\pi)<+\infty$ one obtains
\[
 F'(\pi)\leq\liminf_{j\to\infty}F'(Q_j)
 \leq\liminf_{j\to\infty}F_n^0(Q_j)=F_n^0(\pi),
 \qquad Q_j\in\mathcal D_n, Q_j\weakto\pi.
\]
This is the $\Gamma$-limsup inequality.

It remains to check equicoercivity.  Every sublevel set of $F_n^\varepsilon$, uniformly for
$0<\varepsilon\leq1$, is contained in $\Pi_{A,\mathrm{sym}}^n$.  Since $\Theta$ is coercive,
$\int_Z\Theta\d\pi\leq A$ implies uniform tightness.  Prokhorov's theorem therefore
gives equicoercivity in the narrow topology.

The convergence of minima and the assertion on minimizers follow from the fundamental
theorem of $\Gamma$-convergence (see, for example, \cite[Chapters~7--9]{DalM1989Birk}).
In the present setting one can also see the value
convergence directly: the liminf inequality gives the lower bound, while a recovery
sequence at any minimizer of $\C_{\mathrm{sym}}^n$ gives the upper bound.
\end{proof}

\begin{rmk}[Recovery from the diagonal estimate]
\label{subsec:diagonal-recovery}
For the classical problem, in which the density is prescribed, Lewin \cite{Lewi2018CRM} proved the convergence of the minimal values in the semiclassical limit of the Levy--Lieb functional by constructing a recovery sequence at a minimizer only: an optimal transport plan is supported away from the diagonal, it is regularized at fixed marginals as in \cite{BindDePa2017JEP}, and a fermionic trial state is built from Slater determinants of well-separated orbitals. One may wonder whether the same strategy applies to the moment-constrained problem. In this setting, \Cref{thm:diagonal-estimate} plays the role of the diagonal exclusion for optimal plans.

We explain how \Cref{thm:diagonal-estimate} supplies the recovery sequence needed
for the value convergence.  In the marginal-moment setting of \Cref{sec:repulsive-costs}, with $X=\R^3$ and the Coulomb cost, if $\pi_n^\star$ is a minimizer of the $n$-moment problem and
$\kappa_\rho(\alpha)<1/(2(M-1))$, \Cref{thm:diagonal-estimate} gives
\[
 \pi_n^\star(D_{\alpha_n})=0,
 \qquad
 h(\alpha_n)=
 \frac{I(\rho)-I^n(\rho)}{1-2(M-1)\kappa_\rho(\alpha)}+2(M-1)h(\alpha).
\]
Thus the minimizer has a positive distance from the repulsive diagonal, which, by \Cref{cor:diagonal-uniform}, is bounded from below independently of
the number of moments.  After approximating $\pi_n^\star$ by the finitely supported measures of \Cref{prop:finite-cost-separation}, whose atoms lie in $\supp(\pi_n^\star)$, one may smooth the atoms by packets of radius
$\eta<\alpha_n/3$.  The packets then remain away from the diagonal, and the Coulomb cost is
continuous on their common support.  The usual moment correction restores the constraints
exactly.

More precisely, if this construction gives states $\Gamma_{n,\eta}$ with
\[
 K(\Gamma_{n,\eta})\leq A_n\eta^{-2},
 \qquad
 \pi_{\Gamma_{n,\eta}}\rightharpoonup\pi_n^\star,
 \qquad
 \int_Z c\,\d\pi_{\Gamma_{n,\eta}}\longrightarrow I^n(\rho),
\]
then for fixed $n$ one may take $\eta=\varepsilon^{1/4}$ for all sufficiently small
$\varepsilon$.  This gives $\varepsilon K(\Gamma_{n,\eta})=O(\sqrt\varepsilon)$ and hence
\[
 \limsup_{\varepsilon\downarrow0}\mathcal L_n^\varepsilon
 \leq I^n(\rho).
\]
Together with the lower bound, this proves convergence of the minimum values and of
minimizers.  Notice that this argument uses recovery at minimizers only.  It is therefore
enough for the variational conclusion, but it does not replace the density construction
above when the full $\Gamma$-limsup is asserted.

Finally, since $\alpha_n\geq\alpha^\star$ for every $n$ by \Cref{cor:diagonal-uniform}, the separation part of the construction is
uniform in $n$.  A joint limit $n=n(\varepsilon)\to\infty$ would additionally
require uniform bounds for the moment-correction constants; the diagonal estimate alone does
not provide such bounds.
\end{rmk}

\endgroup

\section{Other applications}\label{app:other-applications}

We now show how the general framework of \Cref{sec:moment-theory} applies to other transport problems: moment inequalities describe partial transport and capacity constraints, additional equalities encode the martingale condition, and the underlying space may itself be a space of probability measures. We state the theorems without proof because their proofs are analogous to those in \Cref{sec:mk}. We highlight any differences in the arguments when necessary.

\subsection{Partial optimal transport}\label{sec:partial-transport}

Let $X$, $Y$ and $c$ be as above, let $\mu\in\mathcal M_+(X)$ and $\nu\in\mathcal M_+(Y)$, and set $m=\min\{\mu(X),\nu(Y)\}$. The partial optimal transport problem with transported mass $m$ is
\begin{equation}\label{POT}
    \inf_{\pi\in\Pi_{\leq}(\mu,\nu)}
    \int_{X\times Y}c(x,y)\,\d\pi(x,y),
\end{equation}
where
\[
\begin{aligned}
    \Pi_{\leq}(\mu,\nu)\coloneqq\{\, &
    \pi\in\mathcal M_+(X\times Y):\ \pi(X\times Y)=m,\\
    &(p_X)_\#\pi\leq\mu,\quad (p_Y)_\#\pi\leq\nu\,\}.
\end{aligned}
\]
If $m=0$, the only admissible measure is the zero measure, and the optimal value is zero. Henceforth assume $m>0$.

To obtain a finite moment relaxation, take $Z=X\times Y$, choose nonnegative functions $\varphi_1,\dots,\varphi_r\in C_b(X)$ and $\psi_1,\dots,\psi_s\in C_b(Y)$, and set $n=r+s$ and
\[
\begin{aligned}
    w_i
    &=\frac1m\int_X\varphi_i\,\d\mu-\varphi_i\circ p_X,
    &&i=1,\dots,r,\\
    w_{r+j}
    &=\frac1m\int_Y\psi_j\,\d\nu-\psi_j\circ p_Y,
    &&j=1,\dots,s.
\end{aligned}
\]
Retain a coercive moment bound to ensure narrow compactness of the finite moment class. Choose $\theta_X$ and $\theta_Y$ with the properties described above and set
\[
    \Theta(x,y)=\theta_X(x)+\theta_Y(y),
    \qquad
    A>0,\quad
    A\geq\int_X\theta_X\,\d\mu+\int_Y\theta_Y\,\d\nu.
\]
Marginal domination then gives
\[
    \int_{X\times Y}\Theta\,\d\pi
    \leq\int_X\theta_X\,\d\mu+\int_Y\theta_Y\,\d\nu
    \leq A
    \qquad\text{for every }\pi\in\Pi_{\leq}(\mu,\nu).
\]
Thus $\Pi_{\leq}(\mu,\nu)\subseteq\Pi_A^n$, and \eqref{eq:moment-problem} provides a relaxation of \eqref{POT}.

With these choices, set $\Pi^n_{A,\leq}(\mu,\nu):=\Pi_A^n$ and consider the finite moment approximation
\begin{equation}\label{MCPOT}
    \inf_{\pi\in\Pi^n_{A,\leq}(\mu,\nu)}
    \int_{X\times Y}c(x,y)\,\d\pi(x,y).
\end{equation}

\begin{thr}\label{thr:existence-MCPOT}
If \eqref{MCPOT} is finite, then it admits a minimizer which is an extreme point of $\Pi^n_{A,\leq}(\mu,\nu)$ and it is supported on at most $r+s+2$ points.
\end{thr}

For the duality result below we introduce the following notation $E_r=\operatorname{cone}_+(\{\varphi_1,\dots,\varphi_r\}),$ and $E_s=\operatorname{cone}_+(\{\psi_1,\dots,\psi_s\}).$ In particular, every function in $-E_r$ or $-E_s$ is nonpositive.

\begin{thr}\label{thr:duality-MCPOT}
Assume that $\Theta$ is finite-valued and continuous. If there exists $\eta\in\Pi^n_{A,\leq}(\mu,\nu)$ such that $\int_{X\times Y}c\,\d\eta<+\infty$ and $\supp(\eta)=\overline{\{c<+\infty\}},$ then
\begin{equation}\label{duality-MCPOT}
\begin{aligned}
    &\min_{\pi\in\Pi^n_{A,\leq}(\mu,\nu)}
    \int_{X\times Y}c\,\d\pi=
    \\&
    \max_{\substack{
        \varphi\in-E_r,\ \psi\in-E_s\\
        \lambda\in\R,\ \tau\geq0\\
        \varphi(x)+\psi(y)
        \leq c(x,y)-\lambda+\tau\Theta(x,y)
    }}
    \left\{
        \int_X\varphi\,\d\mu
        +\int_Y\psi\,\d\nu
        +m\lambda-\tau A
    \right\}.
\end{aligned}
\end{equation}
\end{thr}

\begin{df}\label{df:domination-characterization}
Let $\mu\in\M_+(X)$ and $m>0$. A family of nonnegative functions $(\varphi_i)_{i\in\N}\subset C_b(X)$ characterizes domination by $\mu$ among measures of mass $m$ if, for every $\eta\in\M_+(X)$ with $\eta(X)=m$,
\[
    \left(
    \int_X\varphi_i\,\d\eta\leq\int_X\varphi_i\,\d\mu
    \text{ for every }i\in\N
    \right)
    \Longrightarrow \eta\leq\mu.
\]
\end{df}

\begin{thr}\label{thr:convergence-MCPOT}
Let $(\varphi_i)_{i\in\N}\subset C_b(X)$ and $(\psi_j)_{j\in\N}\subset C_b(Y)$ characterize domination by $\mu$ and $\nu$, respectively, among measures of mass $m$. With $\Theta$ and $A$ chosen as above, let $I^n$ be the minimum in \eqref{MCPOT} using the first $n$ tests of each family, and let $I$ be the minimum in \eqref{POT}. If $I<+\infty$, then $I^n\nearrow I$ as $n\to+\infty$. Moreover, every sequence of minimizers of \eqref{MCPOT} admits a narrowly convergent subsequence, and every accumulation point is a minimizer of \eqref{POT}.
\end{thr}

\begin{cor}\label{cor:duality-POT}
Under the assumptions of \Cref{thr:convergence-MCPOT}, assume that \eqref{POT} is finite and that \Cref{thr:duality-MCPOT} applies for every $n\in\N$. Then
\[
    \min_{\pi\in\Pi_{\leq}(\mu,\nu)}
    \int_{X\times Y}c\,\d\pi
    =
    \sup_{\substack{
        \varphi\in C_b(X),\ \psi\in C_b(Y)\\
        \varphi\leq0,\ \psi\leq0,\ \lambda\in\R\\
        \varphi(x)+\psi(y)\leq c(x,y)-\lambda
    }}
    \left\{
        \int_X\varphi\,\d\mu
        +\int_Y\psi\,\d\nu+m\lambda
    \right\}.
\]
\end{cor}

The dual problem in \Cref{cor:duality-POT} has the same structure as that in \cite[Theorem 2.2]{IN18}.

\subsection{The JKO scheme}\label{sec:jko}

Let $F:\P_2(\R^d)\to\R\cup\{+\infty\}$ be a functional, $\tau>0$, and $\nu\in\P_2(\R^d)$. The associated JKO step is
\begin{equation}\label{JKO}
    \inf_{\rho\in\P_2(\R^d)}
    \left\{
    F(\rho)+\frac{1}{2\tau}W_2^2(\rho,\nu)
    \right\}.
\end{equation}
Assume that $F$ is linear and of the form
\[
    F(\rho)=\int_{\R^d}V(x)\,\d\rho(x),
\]
where $V:\R^d\to\R\cup\{+\infty\}$ is lower semicontinuous, bounded below, and not identically $+\infty$. Define $c_\tau(x,y)=V(x)+\frac{1}{2\tau}|x-y|^2.$ Then \eqref{JKO} is equivalent to
\[
    \inf_{\substack{
        \pi\in\P_2(\R^d\times\R^d)\\
        (p_Y)_\#\pi=\nu
    }}
    \int_{\R^d\times\R^d}c_\tau(x,y)\,\d\pi(x,y),
\]
where $p_Y(x,y)=y$. Thus, the first marginal is free, while the second marginal is prescribed. To obtain a finite moment relaxation, take $Z=\R^d\times\R^d$ and $m=1$, choose $\psi_1,\dots,\psi_s\in C_b(\R^d)$, and define
\[
    g_j=\psi_j\circ p_Y-\int_{\R^d}\psi_j\,\d\nu,
    \qquad j=1,\dots,s.
\]
Using $g_j,-g_j$ as the constraint functions gives $n=2s$ inequalities of the form appearing in \eqref{eq:moment-problem}.

A compatible coercive moment bound can be obtained by taking
\[
    \Theta(x,y)=|x-y|^2+|y|^2.
\]
Indeed, let $b\in\R$ satisfy $V\geq b$, choose $x_0\in\R^d$ with $V(x_0)<+\infty$, and set
\[
    B=V(x_0)+\frac{1}{2\tau}
    \int_{\R^d}|x_0-y|^2\,\d\nu(y).
\]
The plan $\delta_{x_0}\otimes\nu$ has cost $B$. Every plan with second marginal $\nu$ and cost at most $B$ satisfies
\[
    \int_{\R^d\times\R^d}\Theta\,\d\pi
    \leq 2\tau(B-b)+\int_{\R^d}|y|^2\,\d\nu(y).
\]
Consequently, choosing $A>0$ at least as large as the right-hand side preserves the optimal value of \eqref{JKO}. With this choice, \eqref{eq:moment-problem} is a finite moment relaxation of the JKO step, with a narrowly compact admissible set.

Define
\begin{equation}\label{MCJKO-set}
\begin{aligned}
    \Pi_A^n(\nu)\coloneq{}&
    \bigl\{\pi\in\P(\R^d\times\R^d):
    \int_{\R^d\times\R^d}\Theta\,\d\pi\leq A,\\
    &\quad
    \int_{\R^d\times\R^d}\psi_j(y)\,\d\pi(x,y)
    =\int_{\R^d}\psi_j\,\d\nu,
    \quad j=1,\dots,s
    \bigr\}.
\end{aligned}
\end{equation}
We consider the finite moment problem
\begin{equation}\label{MCJKO}
    \inf_{\pi\in\Pi_A^n(\nu)}
    \int_{\R^d\times\R^d}c_\tau(x,y)\,\d\pi(x,y).
\end{equation}

\begin{thr}\label{thr:existence-MCJKO}
If \eqref{MCJKO} is finite, then it admits a minimizer which is an extreme point of $\Pi_A^n(\nu)$ and it is supported on at most $s+2$ points.
\end{thr}

For the duality result below we use $\beta$ for the multiplier associated with the coercive moment bound, since $\tau$ is classically used to denote the time step. Define $E_s =\span(\{1,\psi_1,\dots,\psi_s\})$.

\begin{thr}\label{thr:duality-MCJKO}
Assume that $\Theta$ is finite-valued and continuous. If there exists $\eta\in\Pi_A^n(\nu)$ such that $\int_{\R^d\times\R^d}c_\tau\,\d\eta<+\infty$ and $\supp(\eta)=\overline{\{c_\tau<+\infty\}},$ then
\begin{equation}\label{duality-MCJKO}
\begin{aligned}
    \min_{\pi\in\Pi_A^n(\nu)}
    \int_{\R^d\times\R^d}c_\tau\,\d\pi
    =
    \max_{\substack{
        \psi\in E_s,\ \beta\geq0\\
        \psi(y)\leq c_\tau(x,y)+\beta\Theta(x,y)
    }}
    \left\{
        \int_{\R^d}\psi\,\d\nu-\beta A
    \right\}.
\end{aligned}
\end{equation}
\end{thr}

\begin{thr}\label{thr:convergence-MCJKO}
Let $(\psi_j)_{j\in\N}\subset C_b(\R^d)$ fully characterize $\nu$ in the sense of \Cref{def:moment-determinacy}. Let $I^n$ be the minimum in \eqref{MCJKO} using the first $n$ tests, and let $I$ be the minimum in \eqref{JKO}. Assume that $I<+\infty$ and that some optimal plan for the lifted formulation of \eqref{JKO} satisfies $\int\Theta\,\d\pi\leq A$. Then $I^n\nearrow I$ as $n\to+\infty$. Moreover, every sequence of minimizers of \eqref{MCJKO} admits a narrowly convergent subsequence, and every accumulation point $\pi$ satisfies $(p_Y)_\#\pi=\nu$ and has first marginal $(p_X)_\#\pi$ minimizing \eqref{JKO}.
\end{thr}

\begin{cor}\label{cor:duality-JKO}
Under the assumptions of \Cref{thr:convergence-MCJKO}, assume that \Cref{thr:duality-MCJKO} applies for every $n\in\N$. Then
\[
    \min_{\rho\in\P_2(\R^d)}
    \left\{
        \int_{\R^d}V\,\d\rho
        +\frac{1}{2\tau}W_2^2(\rho,\nu)
    \right\}
    =
    \sup_{\substack{
        \psi\in C_b(\R^d), \beta \geq 0\\
        \psi(y)\leq c_\tau(x,y) + \beta \Theta(x,y)} }
    \int_{\R^d}\psi\,\d\nu - \beta A.
\]
\end{cor}

Our dual formula differs in appearance from \cite[equation (1.4)]{JLL21}, where the maximization involves two potentials and the convex conjugate of the energy. Here, only one potential appears, and the energy enters directly through the pointwise constraint
\[
    \psi(y)\leq V(x)+\frac{1}{2\tau}|x-y|^2.
\]
This simplification results from the fact that the contribution of a linear energy can be incorporated into the transport cost, allowing the auxiliary potential in their formulation to be eliminated.

\subsection{Martingale optimal transport}\label{sec:martingale-transport}

Let $\mu,\nu\in\P_1(\R^d)$, where $\P_1(\R^d)$ denotes the probability measures with finite first moment, and let $c:\R^d\times\R^d\to\R\cup\{+\infty\}$ be lower semicontinuous and bounded below. For $\pi\in\Pi(\mu,\nu)$, write $\pi(\d x,\d y)=\mu(\d x)\pi_x(\d y)$ for its disintegration with respect to the first marginal. The martingale optimal transport problem is
\begin{equation}\label{martOT}
    \inf_{\substack{
        \pi\in\Pi(\mu,\nu)\\
        \int_{\R^d}y\,\d\pi_x(y)=x\ \mu\text{-a.e.}
    }}
    \int_{\R^d\times\R^d}c(x,y)\,\d\pi(x,y).
\end{equation}
Its admissible set is nonempty if and only if $\mu\preceq_{\mathrm{cx}}\nu$, where $\preceq_{\mathrm{cx}}$ denotes convex order (see \cite{AlfCoyEhrLom2021Approx} for details).

Take $Z=\R^d\times\R^d$ and $m=1$. The martingale condition is equivalent to
\[
    \int_{\R^d\times\R^d}
    \phi(x)\cdot(y-x)\,\d\pi(x,y)=0
    \qquad
    \text{for every }\phi\in C_b(\R^d;\R^d).
\]
Thus, in addition to the marginal moment constraints, one may select finitely many functions $\phi_1,\dots,\phi_q$ and impose the inequalities associated with $h_\ell$ and $-h_\ell$, where
\[
    h_\ell(x,y)=\phi_\ell(x)\cdot(y-x),
    \qquad \ell=1,\dots,q.
\]
These functions are continuous but generally unbounded. Consequently, this falls outside the assumption $w_i\in C_b(Z)$ and requires an extension to continuous unbounded tests, together with suitable uniform integrability conditions.

For exact marginal constraints, such an extension is considered in \cite{Zae2015Monge}. For Polish spaces $X$ and $Y$, the author introduces the space $C_L(\mu,\nu)$ of functions $h\in C(X\times Y)$ satisfying $|h(x,y)|\leq f(x)+g(y)$ for some $f\in C(X)\cap L^1(\mu)$ and $g\in C(Y)\cap L^1(\nu)$, equipped with the seminorm
\[
    \|h\|_L
    \coloneqq
    \sup_{\pi\in\Pi(\mu,\nu)}
    \int_{X\times Y}|h|\,\d\pi.
\]
The space $C_b(X\times Y)$ is dense in $C_L(\mu,\nu)$ for this seminorm, and $\pi\mapsto\int h\,\d\pi$ is narrowly continuous on $\Pi(\mu,\nu)$ for every $h\in C_L(\mu,\nu)$. In particular, the martingale tests belong to this space, since
\[
    |\phi(x)\cdot(y-x)|
    \leq \|\phi\|_\infty(|x|+|y|).
\]
This continuity statement relies on the exact marginal constraints and does not automatically extend to their finite moment relaxation.

For the latter, one may impose growth conditions of the type used in \cite{AlfCoyEhrLom2021Approx}. Specifically, for the selected marginal tests $\varphi_i,\psi_j$ and martingale tests $h_\ell$, a sufficient condition is
\[
    |\varphi_i(x)|+|\psi_j(y)|+|h_\ell(x,y)|
    \leq C\bigl(1+\Theta(x,y)\bigr)^s,
    \qquad 0<s<1,
\]
for some $C>0$, where
\[
    \Theta(x,y)=\theta_\mu(|x|)+\theta_\nu(|y|).
\]
Here $\theta_\mu,\theta_\nu:\R_+\to\R_+$ are continuous, nondecreasing functions tending to infinity, with
\[
    \int_{\R^d}\theta_\mu(|x|)\,\d\mu(x)
    +\int_{\R^d}\theta_\nu(|y|)\,\d\nu(y)<\infty.
\]
These are additional assumptions on the selected tests and the marginals. Choose $A>0$ at least as large as this sum. The bound $\int\Theta\,\d\pi\leq A$, together with the sublinear growth condition, ensures uniform integrability of the constraint functions. Their moment constraints therefore pass to narrow limits, allowing the finite moment formulation to accommodate these unbounded tests.

For the next results, assume that $X,Y\subset\R^d$ are compact, so that we may take $\Theta\equiv0$. In addition to the finite marginal constraints, choose $H_1,\dots,H_q\in C(X;\R^d)$ and impose
\[
    \int_{X\times Y}H_\ell(x)\cdot(y-x)\,\d\pi(x,y)=0,
    \qquad \ell=1,\dots,q.
\]
The functions $h_\ell(x,y)=H_\ell(x)\cdot(y-x)$ belong to $C_b(X\times Y)$, so these equalities are encoded by adding $h_\ell,-h_\ell$ to the moment functions. We denote the resulting admissible set by $\Pi^n_{\mathrm{mart}}(\mu,\nu)$, with $n=2(r+s+q)$, and consider
\begin{equation}\label{MCMOT}
    \inf_{\pi\in\Pi^n_{\mathrm{mart}}(\mu,\nu)}
    \int_{X\times Y}c(x,y)\,\d\pi(x,y).
\end{equation}

\begin{thr}\label{thr:existence-MCMOT}
If \eqref{MCMOT} is finite, then it admits a minimizer which is an extreme point of $\Pi^n_{\mathrm{mart}}(\mu,\nu)$ and is supported on at most $r+s+q+1$ points.
\end{thr}

For the duality result below we introduce the following notation
\[
    E_r=\operatorname{span}\{1,\varphi_1,\dots,\varphi_r\},
    \qquad
    E_s=\operatorname{span}\{1,\psi_1,\dots,\psi_s\},
    \qquad
    \mathcal H_q=\operatorname{span}\{H_1,\dots,H_q\}.
\]

\begin{thr}\label{thr:duality-MCMOT}
If there exists $\eta\in\Pi^n_{\mathrm{mart}}(\mu,\nu)$ such that $\int_{X\times Y}c\,\d\eta<+\infty$ and $\supp(\eta)=\overline{\{c<+\infty\}},$ then
\begin{equation}\label{duality-MCMOT}
\begin{aligned}
    &\min_{\pi\in\Pi^n_{\mathrm{mart}}(\mu,\nu)}
    \int_{X\times Y}c\,\d\pi
    \quad=
    \max_{\substack{
        \varphi\in E_r,\ \psi\in E_s,\ H\in\mathcal H_q\\
        \varphi(x)+\psi(y)+H(x)\cdot(y-x)\leq c(x,y)
    }}
    \left\{
        \int_X\varphi\,\d\mu+\int_Y\psi\,\d\nu
    \right\}.
\end{aligned}
\end{equation}
\end{thr}

\begin{thr}\label{thr:convergence-MCMOT}
Let $(\varphi_i)_{i\in\N}\subset C(X)$ and $(\psi_j)_{j\in\N}\subset C(Y)$ fully characterize $\mu$ and $\nu$, respectively, and let $(H_\ell)_{\ell\in\N}\subset C(X;\R^d)$ have a uniformly dense linear span in $C(X;\R^d)$. Let $I^k$ be the minimum in \eqref{MCMOT} using the first $k$ tests of each family, so that $r=s=q=k$ and $n=6k$, and let $I$ be the minimum in \eqref{martOT}. If $I<+\infty$, then $I^k\nearrow I$ as $k\to+\infty$. Moreover, every sequence of minimizers of \eqref{MCMOT} admits a narrowly convergent subsequence, and every accumulation point is a minimizer of \eqref{martOT}.
\end{thr}

\begin{cor}\label{cor:duality-MOT}
Under the assumptions of \Cref{thr:convergence-MCMOT}, assume that \eqref{martOT} is finite and that \Cref{thr:duality-MCMOT} applies at every approximation level. Then
\[
    \min_{\pi\in\Pi_{\mathrm{mart}}(\mu,\nu)}
    \int_{X\times Y}c\,\d\pi
    =
    \sup_{\substack{
        \varphi\in C(X),\ \psi\in C(Y)\\
        H\in C(X;\R^d)\\
        \varphi(x)+\psi(y)+H(x)\cdot(y-x)\leq c(x,y)
    }}
    \left\{
        \int_X\varphi\,\d\mu+\int_Y\psi\,\d\nu
    \right\},
\]
where $\Pi_{\mathrm{mart}}(\mu,\nu)$ denotes the admissible set in \eqref{martOT}.
\end{cor}

The dual formulation above has the same structure as in \cite{BHLP13,GO19}.

\subsection{Capacity constrained optimal transport}\label{sec:capacity-transport}

Let $X$, $Y$ and $c$ be as in the Monge--Kantorovich problem, let $\mu\in\P(X)$, $\nu\in\P(Y)$, and let $\sigma\in\mathcal M_+(X\times Y)$ be a prescribed capacity measure. Consider
\begin{equation}\label{CCOT}
    \inf_{\substack{\pi\in\Pi(\mu,\nu)\\ \pi\leq\sigma}}
    \int_{X\times Y}c(x,y)\,\d\pi(x,y).
\end{equation}
This includes a prescribed upper bound on the transport density by taking $\sigma=\bar h\,\lambda$ for a reference measure $\lambda$ and a nonnegative $\bar h\in L^1(\lambda)$. The Euclidean density formulation is studied in \cite{KM15}.

Take $Z=X\times Y$ and $m=1$. The capacity constraint is equivalent to
\[
    \int_{X\times Y}u\,\d\pi
    \leq\int_{X\times Y}u\,\d\sigma
    \qquad
    \text{for every }u\in C_b(X\times Y),\quad u\geq0.
\]
Thus, choose finitely many nonnegative functions $u_1,\dots,u_q\in C_b(X\times Y)$ and define
\[
    a_k=\int_{X\times Y}u_k\,\d\sigma-u_k,
    \qquad k=1,\dots,q.
\]
Adding the inequalities $\int a_k\,\d\pi\geq0$ to the finite marginal constraints gives a problem of the form \eqref{eq:moment-problem}. Retaining the compatible coercive moment bound constructed from $\mu$ and $\nu$ ensures narrow compactness and preserves every admissible plan of \eqref{CCOT}.

Define, for $n=2(r+s)+q$,
\begin{equation}\label{MCCOT-set}
\begin{aligned}
    \Pi^n_{A,\sigma}(\mu,\nu)\coloneq{}&
    \bigl\{\pi\in\Pi_A^{2(r+s)}(\mu,\nu): \quad
    \int_{X\times Y}a_k\,\d\pi\geq0,
    \quad k=1,\dots,q
    \bigr\},
\end{aligned}
\end{equation}
where $\Pi_A^{2(r+s)}(\mu,\nu)$ is the marginal moment class introduced in the Monge--Kantorovich case. We consider the finite moment approximation of \eqref{CCOT}
\begin{equation}\label{MCCOT}
    \inf_{\pi\in\Pi^n_{A,\sigma}(\mu,\nu)}
    \int_{X\times Y}c(x,y)\,\d\pi(x,y).
\end{equation}

\begin{thr}\label{thr:existence-MCCOT}
If \eqref{MCCOT} is finite, then it admits a minimizer which is an extreme point of $\Pi^n_{A,\sigma}(\mu,\nu)$ and is supported on at most $r+s+q+2$ points.
\end{thr}

For the duality result below we consider $E_r$ and $E_s$ as in the Monge--Kantorovich case and we introduce also $U_q=\operatorname{cone}_+(\{u_1,\dots,u_q\}).$

\begin{thr}\label{thr:duality-MCCOT}
Assume that $\Theta$ is finite-valued and continuous. If there exists $\eta\in\Pi^n_{A,\sigma}(\mu,\nu)$ such that $\int_{X\times Y}c\,\d\eta<+\infty$ and $\supp(\eta)=\overline{\{c<+\infty\}},$ then
\begin{equation}\label{duality-MCCOT}
\begin{aligned}
    &\min_{\pi\in\Pi^n_{A,\sigma}(\mu,\nu)}
    \int_{X\times Y}c\,\d\pi =
    \max_{\substack{
        \varphi\in E_r,\ \psi\in E_s\\
        u\in U_q,\ \tau\geq0\\
        \varphi+\psi\leq c+u+\tau\Theta
    }}
    \left\{
        \int_X\varphi\,\d\mu
        +\int_Y\psi\,\d\nu
        -\int_{X\times Y}u\,\d\sigma
        -\tau A
    \right\}.
\end{aligned}
\end{equation}
\end{thr}

\begin{rmk}
The capacity constraint is imposed through finitely many moment inequalities, so measures in $\Pi^n_{A,\sigma}(\mu,\nu)$ need not satisfy $\pi\leq\sigma$. In particular, they may be atomic even if $\sigma$ has no atoms, whereas the constraint $\pi\leq\sigma$ would rule out atoms in this case.
\end{rmk}

\begin{thr}\label{thr:convergence-MCCOT}
Let $(\varphi_i)_{i\in\N}\subset C_b(X)$ and $(\psi_j)_{j\in\N}\subset C_b(Y)$ fully characterize $\mu$ and $\nu$, respectively, in the sense of \Cref{def:moment-determinacy}. Let $(u_k)_{k\in\N}\subset C_b(X\times Y)$ be nonnegative and characterize domination by $\sigma$ among measures of mass $1$, in the sense of \Cref{df:domination-characterization}. With $\Theta$ and $A$ chosen as above, let $I^n$ be the minimum in \eqref{MCCOT} using the first $n$ tests of each family, and let $I$ be the minimum in \eqref{CCOT}. If $I<+\infty$, then $I^n\nearrow I$ as $n\to+\infty$. Moreover, every sequence of minimizers of \eqref{MCCOT} admits a narrowly convergent subsequence, and every accumulation point is a minimizer of \eqref{CCOT}.
\end{thr}

\begin{cor}\label{cor:duality-CCOT}
Under the assumptions of \Cref{thr:convergence-MCCOT}, assume that \Cref{thr:duality-MCCOT} applies for every $n\in\N$. Then
\[
    I=
    \sup_{\substack{\varphi\in C_b(X)\\ \psi\in C_b(Y)}}
    \left\{
        \int_X\varphi\,\d\mu
        +\int_Y\psi\,\d\nu
        -\int_{X\times Y}
        \bigl(\varphi(x)+\psi(y)-c(x,y)\bigr)_+\,\d\sigma
    \right\}.
\]
Here $a_+=\max\{a,0\}$, with the convention $(-\infty)_+=0$.
\end{cor}

Our dual formula is similar to \cite[eqs.~(4)--(6)]{KMS15} with $\sigma = \bar h \d x \d y$ and surplus $s=-c$. Therefore, our expression is the negative of theirs.

\subsection{Wasserstein-on-Wasserstein transport}\label{sec:wasserstein-on-wasserstein}

Let $(S,d)$ be a complete separable metric space, $1\leq p<\infty$, and set $E=(\P_p(S),W_p)$, which is itself a Polish space. For two laws $\mathbb M,\mathbb N\in\P_p(E)$, the Wasserstein-on-Wasserstein distance is defined by
\begin{equation}\label{WoW}
    \mathbb W_p^p(\mathbb M,\mathbb N) = \inf_{\Gamma\in\Pi(\mathbb M,\mathbb N)} \int_{E\times E}W_p^p(\rho,\sigma)\,\d\Gamma(\rho,\sigma).
\end{equation}
See \cite{PinzSava2025arXiv}. Here the transported objects are probability measures, and the objective remains linear in the coupling $\Gamma$.

This problem fits our framework with $Z=E\times E$, $m=1$ and $c(\rho,\sigma)=W_p^p(\rho,\sigma).$ Choose $F_1,\dots,F_r,G_1,\dots,G_s\in C_b(E)$ and define
\[
\begin{aligned}
    f_i&=F_i\circ p_1-\int_EF_i\,\d\mathbb M,
    &&i=1,\dots,r,\\
    g_j&=G_j\circ p_2-\int_EG_j\,\d\mathbb N,
    &&j=1,\dots,s,
\end{aligned}
\]
where $p_1,p_2:E\times E\to E$ are the canonical projections. Using $f_i,-f_i,g_j,-g_j$ as the constraint functions gives the finite moment relaxation.

As before, choose nonnegative lower semicontinuous functions $\theta_{\mathbb M},\theta_{\mathbb N}:E\to[0,+\infty]$ with compact sublevel sets in the $W_p$ topology and finite integrals against the respective laws. Set
\[
\begin{aligned}
    \Theta(\rho,\sigma)
    &=\theta_{\mathbb M}(\rho)+\theta_{\mathbb N}(\sigma),\\
    A&\geq
    \int_E\theta_{\mathbb M}\,\d\mathbb M
    +\int_E\theta_{\mathbb N}\,\d\mathbb N,
    \qquad A>0.
\end{aligned}
\]
This bound preserves every exact coupling and ensures narrow compactness of the finite moment class.

The tests $F_i$ and $G_j$ are functionals on the space of measures. Restricting them to linear evaluations $\rho\mapsto\int_S\varphi\,\d\rho$ would constrain only the barycenters of the marginal laws, rather than determine the laws themselves.

We denote by $\Pi_A^n(\mathbb M,\mathbb N)$ the admissible set corresponding to the chosen marginal tests, with $n=2(r+s)$, and consider
\begin{equation}\label{MCWoW}
    \inf_{\Gamma\in\Pi_A^n(\mathbb M,\mathbb N)}
    \int_{E\times E}W_p^p(\rho,\sigma)\,\d\Gamma(\rho,\sigma).
\end{equation}

\begin{thr}\label{thr:existence-MCWoW}
Problem \eqref{MCWoW} admits a minimizer which is an extreme point of $\Pi_A^n(\mathbb M,\mathbb N)$ and is supported on at most $r+s+2$ points of $E\times E$.
\end{thr}

For the duality result let $E_r=\span(\{F_0,\dots,F_r\})$ and $E_s=\span(\{G_0,\dots,G_s\})$ where $F_0\equiv G_0\equiv1$.

\begin{thr}\label{thr:duality-MCWoW}
Assume that $S$ is compact and choose $\Theta\equiv0$. If there exists $\eta\in\Pi_A^n(\mathbb M,\mathbb N)$ such that $\supp(\eta)=E\times E,$ then
\begin{equation}\label{duality-MCWoW}
\begin{aligned}
    &\min_{\Gamma\in\Pi_A^n(\mathbb M,\mathbb N)}
    \int_{E\times E}W_p^p(\rho,\sigma)\,\d\Gamma(\rho,\sigma)=
    \max_{\substack{
        F\in E_r,\ G\in E_s\\
        F(\rho)+G(\sigma)\leq W_p^p(\rho,\sigma)
    }}
    \left\{
        \int_E F\,\d\mathbb M
        +\int_E G\,\d\mathbb N
    \right\}.
\end{aligned}
\end{equation}
\end{thr}

\begin{rmk}
The discrete characterization concerns the coupling $\Gamma$ on $E\times E$: the measures $\rho$ and $\sigma$ appearing in its support need not themselves be discrete.
\end{rmk}

\begin{rmk}
Since $S$ is compact, $E$ is compact and the cost $W_p^p$ is bounded, so finiteness of the cost of $\eta$ is automatic. If $\mathbb M$ and $\mathbb N$ have full support on $E$, one may choose $\eta=\mathbb M\otimes\mathbb N$. When $S$ is noncompact, a finite-valued continuous coercive function $\Theta$ on $E\times E$ need not exist (\cite[p. 978]{LV09}). Indeed, the existence of such a function would imply that $E\times E$, and hence $E$, is locally compact.
\end{rmk}

\begin{thr}\label{thr:convergence-MCWoW}
Let $(F_i)_{i\in\N},(G_j)_{j\in\N}\subset C_b(E)$ fully characterize $\mathbb M$ and $\mathbb N$, respectively, in the sense of \Cref{def:moment-determinacy}. With $\Theta$ and $A$ chosen as above, let $I^n$ be the minimum in \eqref{MCWoW} using the first $n$ tests of each family. Then $I^n\nearrow\mathbb W_p^p(\mathbb M,\mathbb N)\text{ as }n\to+\infty.$ Moreover, every sequence of minimizers of \eqref{MCWoW} admits a narrowly convergent subsequence in $\P(E\times E)$, and every accumulation point is a minimizer of \eqref{WoW}.
\end{thr}

\begin{cor}\label{cor:duality-WoW}
Under the assumptions of \Cref{thr:convergence-MCWoW}, assume that $S$ is compact and take $\Theta\equiv0$. If \Cref{thr:duality-MCWoW} applies for every $n\in\N$, then
\[
    \mathbb W_p^p(\mathbb M,\mathbb N)
    =
    \sup_{\substack{
        F,G\in C_b(E)\\
        F(\rho)+G(\sigma)\leq W_p^p(\rho,\sigma)
        \text{ for every }\rho,\sigma\in E
    }}
    \left\{
        \int_E F\,\d\mathbb M
        +\int_E G\,\d\mathbb N
    \right\}.
\]
\end{cor}

The dual formula above coincides with Kantorovich duality applied to the metric space $(\P_p(S),W_p)$; compare \cite[equation (2.7)]{PinzSava2025arXiv}.

\section*{A.I. Disclosure} In this paper both Claude and ChatGPT were used for grammar and spelling corrections, as well as a first check for the correctness of the proofs. \Cref{thm:wp-winf-bound} was originally obtained with a worst constant for $h_n$. Discussions with ChatGPT helped to improve the constant, hence obtaining the current argument. 
The authors take full responsibility for correctness and the results of this paper.

\section*{Acknowledgements.}
L.N. benefited  from the ANR project GOTA (ANR-23-CE46-0001) and also supported by the French National Research Agency (ANR) under the France 2030 grant reference number ANR-24-RRII-0002 operated by the Inria Quadrant Program (project PiQuanTO).

\bibliographystyle{plain}
\bibliography{MyBiblio-LDP}
\end{document}